\documentclass[12pt]{article}
\usepackage{amsmath,enumerate,amsfonts,color,amssymb,amsthm}

\usepackage{graphicx}

\usepackage{mathptmx}

\usepackage{tikz}

\def\NN{{\mathbb N}}

\def\RR{{\mathbb R}}

\def\Sphere{{\mathbb S}}

\def\mcM{{\mycal M}}
\def\mcF{{\mycal F}}

\def\hh{{\rm \overline H}}
\def\bR{{\mathbb R}}

\def\eps{{\varepsilon}}
\def\loc{{\rm loc}}

\newtheorem{theorem} {\sc  Theorem\rm} [section]

\newtheorem{corollary} [theorem] {\sc  Corollary\rm}
\newtheorem{lemma} [theorem] {\sc  Lemma\rm}
\newtheorem{proposition} [theorem] {\sc  Proposition\rm}

\newtheorem{definition}[theorem]{\sc  Definition\rm}

\newtheorem{remark}{\sc  Remark\rm}[section]

\newtheorem{claim}{Claim}

\def\nd{\noindent}

\newcounter{marnote}

\DeclareFontFamily{OT1}{rsfs}{}
\DeclareFontShape{OT1}{rsfs}{m}{n}{ <-7> rsfs5 <7-10> rsfs7 <10-> rsfs10}{}
\DeclareMathAlphabet{\mycal}{OT1}{rsfs}{m}{n}

\def\tr{{\rm tr}}

\def\mcS{{\mycal{S}}}

\def\hh{{\rm \overline H}}
\def\bR{{\mathbb R}}

\def\be{\begin{equation}}
\def\ee{\end{equation}}

\def\tr{{\rm tr}}

\def\mcS{{\mycal{S}}}

\def\ga{{\gamma_+}}
\def\Spec{\rm Spec}
\def \f {\varphi}
\def\Rot{{\mathcal R}}

\newcommand{\defeq}{\stackrel{\scriptscriptstyle \text{def}}{=}}
\numberwithin{equation}{section}
\begin{document}
\date{}
\title{Uniqueness results for an ODE related to a generalized Ginzburg-Landau
model for liquid crystals}
\author{Radu Ignat\thanks{Institut de Math\'ematiques de Toulouse, Universit\'e
Paul Sabatier, b\^at. 1R3, 118 Route de Narbonne, 31062
Toulouse, France. Email: Radu.Ignat@math.univ-toulouse.fr
}~, Luc Nguyen\thanks{Mathematics Department, Princeton University, Fine Hall, Washington Road, Princeton, NJ 08544, USA. Email: llnguyen@math.princeton.edu}~, Valeriy Slastikov\thanks{School of Mathematics, University of Bristol, Bristol, BS8 1TW, United Kingdom. Email: Valeriy.Slastikov@bristol.ac.uk}~ and Arghir Zarnescu\thanks{University of Sussex, Department of Mathematics, Pevensey 2, Falmer, BN1 9QH, United Kingdom. Email: A.Zarnescu@sussex.ac.uk}}

\maketitle
\begin{abstract}
We study a singular nonlinear ordinary differential equation on intervals $[0,R)$ with $R\le +\infty$, motivated by the Ginzburg-Landau models in superconductivity and Landau-de Gennes models in liquid crystals. We prove existence and uniqueness of positive solutions under general assumptions on the nonlinearity. Further uniqueness results for sign-changing solutions are obtained  for a  physically relevant class of nonlinearities. Moreover, we prove a number of fine qualitative properties of the solution that are important for the study of energetic stability.

\end{abstract}


\section{Introduction}
\label{sect:intropar}

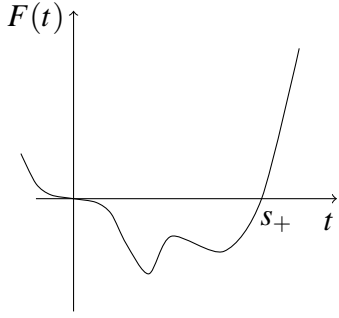
\begin{figure}[h]
\begin{center}
\begin{tikzpicture}
\draw[->] (-0.5,0) -- (3.5,0);
\draw[->] (0,-1.5)--(0,2.5);
\draw plot[smooth] coordinates{(-0.7,0.6)(-0.5,0.2) (-0.3,0.05)   (0,0) (0.3,-0.05) (0.5,-0.2) (0.7,-0.6) (1,-1) (1.3,-0.5) (2,-0.7) (2.5,0) (3,2)};
\draw (3.4,-0.3) node {$t$}
(2.7,-0.3) node {$s_+$}
(-0.5,2.4) node {$F(t)$};
\end{tikzpicture}
\end{center}
\caption{A graph of a prototypical nonlinearity $F$.}
\label{fig:fig1}
\end{figure}

We consider the following ordinary differential equation:
 \begin{align}
&u''(r) + \frac{p}{r}\,u'(r) - \frac{q}{r^2}\,u(r) = F(u(r)) \, \, \text{ in } \, \, (0,R),
\label{RS::ODE}\\
&u(0)=0, \quad u(R) =s_+,\label{BC}
\end{align}  
where $R\le +\infty$, $p$ and $q$ are constants satisfying
\begin{equation}\label{cond:pq}
p,q \in \RR, \ q > 0,
\end{equation}
and $F:\RR\to\RR$ is a $C^1$ function which vanishes at $0$ and at $s_+>0$ (see Fig. \ref{fig:fig1}). In \eqref{BC}, we use the standard convention $u(+\infty):= \lim_{r \to +\infty} u(r) =s_+$ if $R=+\infty$.

The ODE \eqref{RS::ODE} is the Euler-Lagrange equation of the energy functional:
\be\label{energyODE}
E[u;I] = \frac{1}{2}\int_{I}\Big[r^p|u'(r)|^2 + q \, r^{p-2}\,u^2(r) + r^p h(u(r))\Big]\, dr,
\ee
where $I\subset [0, +\infty)$ is an arbitrary interval and
\begin{equation}
h(t):=2\int_0^t F(s)\, ds, \quad t\in \RR.
	\label{Eq:hDefinition}
\end{equation}

The main aim of this paper is to study the existence, uniqueness and qualitative properties of solutions to the boundary value problem \eqref{RS::ODE}, \eqref{BC}. The main difficulty in exploring the ODE satisfied by $u$ is the general type of nonlinearity $F(u)$ on the right hand side. For example, existing techniques for dealing with equations of the type \eqref{RS::ODE} in \cite{tang, Hervex2, Mironescu-radial} are not applicable in our setting. One way of appreciating the effect of the nonlinearity is by noting that for $u\in [0,s_+]$, the function $F$ does not, in general,  satisfy the Krasnosel'ski\v{i} condition (see e.g. \cite{brezis-oswald, krasnoselskii64}), unlike in the standard Ginzburg-Landau case \cite{Mironescu-radial}. Furthermore the Pohozaev-type approach frequently used for proving uniqueness fails in this case.

We start by stating our existence and uniqueness result in the class of non-negative solutions, which was announced in \cite{CRAS-INSZ}.
\begin{theorem}\label{thm:main}
Assume that $p, q$ are given constants satisfying \eqref{cond:pq} and $F: \RR \rightarrow \RR$ is a $C^1$ function satisfying 
\be
\label{condF} 
\begin{cases}
& F(0)=F(s_+)=0, \, \, F'(s_+)>0, \\
& F(t)<0 \, \, \textrm{ if } \, \, t\in (0, s_+), \quad F(t)\ge 0 \, \, \textrm{ if } \, \, t\in(s_+,+\infty).
\end{cases}
\ee
Then there exists a non-negative solution $u$ of the boundary value problem \eqref{RS::ODE} and \eqref{BC}, which is unique in the class of non-negative solutions. Moreover, this solution is strictly increasing.
\end{theorem}

If in addition we restrict the class of nonlinearities, we can show variational properties of the solution:
 \begin{corollary}\label{cor:energyuniqueness}
Assume that $p, q$ are given constants satisfying \eqref{cond:pq} and $F: \RR \rightarrow \RR$ is a $C^1$ function satisfying \eqref{condF} and 
\begin{equation}
F_{even}(t) := \frac{F(t) + F(-t)}{2} \leq 0 \text{ for } t \geq 0.
	\label{Fcond:even}
\end{equation}
Then the solution $u$ in Theorem~\ref{thm:main} is locally energy minimizing with respect to the energy $E$ in \eqref{energyODE}, i.e.
$$E[u;\omega]\le E[u+\varphi;\omega] \, \textrm{ for any } \, \omega\subset [0,+\infty)   \, \textrm{ compact interval and } \, \varphi\in C^\infty_c(\omega).$$ Conversely, if  a function $u\in H^1_{loc}(0,R)$  is locally energy minimizing with respect to $E$ and satisfies $u(R)=s_+$, then $u$ is necessarily the non-negative solution of \eqref{RS::ODE} and \eqref{BC} obtained in Theorem~\ref{thm:main}. 
\end{corollary} 

The proof of Theorem~\ref{thm:main} is split into two parts: existence and uniqueness. The  existence part is done by  constructing energy minimizing solutions on finite intervals and letting the length  of the interval tend to infinity in the case $R=+\infty$. Fine local estimates of the behavior of $u$ near the origin combined with an energy argument ensures the non-flattening of the solution obtained in this limit. The uniqueness part is more delicate to prove. To do this, we construct  comparison barriers through a scaling argument and use suitable versions of the maximum principle together with a detailed understanding of the asymptotics at the origin and at infinity (in the case $R=+\infty$). 
 
One can further ask if the uniqueness result holds for {\it nodal} solutions (i.e. solutions that may change signs). In general, if one assumes only \eqref{condF} then in addition to a  non-negative solution there might exist sign-changing solutions, see Proposition~\ref{lem:MP}. However, under additional assumptions, relevant to the physical problem detailed in subsection~\ref{intro:physical}, we prove
the following uniqueness result for nodal solutions.

\newpage
\begin{figure}[h]
\begin{center}
\includegraphics[scale=0.2]{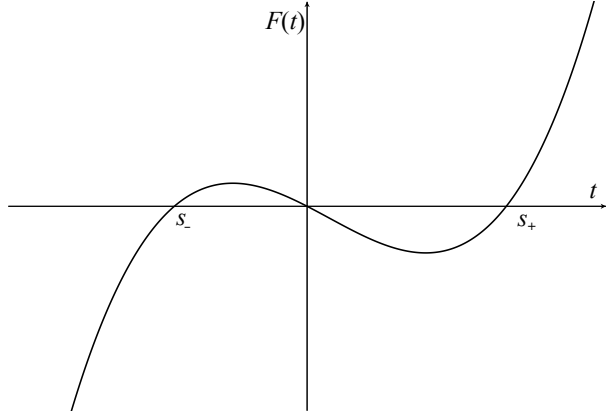}
\end{center}
\caption{A graph of a \emph{physically relevant} nonlinearity $F$.}
\label{fig:fig1+}
\end{figure}

\begin{theorem}\label{thm:uniqueness+}
Assume that $p>0, q>0$ are given constants  and $F: \RR \rightarrow \RR$ is a $C^1$ function satisfying \eqref{condF}. Assume in addition that  there exists   
$s_- \in [-s_+,0)$ such that:
\be
\label{condFLeft}
\begin{cases}
&F(t) \le 0 \, \, \text{ if } \, \, t \in (-\infty,s_-), F(t) \geq 0 \, \, \text{ if }\, \,  t \in (s_-, 0),\vspace{0.2cm}\\
& \displaystyle \frac{F(t_1)}{t_1} + \frac{F(-t_2)}{t_2} \leq 0  \, \, \text{ if } \, \,  0 < t_1 \leq t_2 \leq |s_-|.\\
\end{cases}
\ee
 Then there exists a unique solution $u$ of the boundary value problem \eqref{RS::ODE} and \eqref{BC}.
 \end{theorem}

\begin{remark}
We also prove the above uniqueness result when $p=0$ under additional assumptions on nonlinearity $F$: either we assume in addition that $F$ is a $C^2$ function (see Remark~\ref{rmk:p01}), or
we impose a stronger version of \eqref{condFLeft} for the $C^1$ function $F$, namely,  there exists $\alpha>1$ such that 
\be\label{condF+}
\frac{F(t_1)}{t_1}+\alpha^2 \frac{F(-t_2)}{t_2}\le 0 \  \text{ for every } 0<t_1\le t_2\le |s_-|
\ee
(see Remark~\ref{rmk:p02}). 
For $p < 0$, numerical simulations (see Figure \ref{fig:multiplesol}) suggest that the uniqueness result in Theorem \ref{thm:uniqueness+} does not hold in general (see Remark \ref{Rem:pNeg}).
\end{remark}

\begin{remark}\label{rmk:physrescaling}
The physically relevant  nonlinearity (see Section~\ref{intro:physical}) of the form 
\begin{equation}\label{def:physical nonlinearity}
F(t) =-a^2\,t - \frac{b^2}{3}\,t^2 + \frac{2c^2}{3}\,t^3, \quad t\in \RR
\end{equation} satisfies \eqref{condF} and \eqref{condFLeft} if $a^2, c^2>0$ and $b^2\geq 0$. 
In particular, for $F(t)=-t+t^3$ ($t\in \RR$) and $p=1$ and $q=n^2,n\in\mathbb{Z}\setminus\{0\}$ in \eqref{RS::ODE}, we recover the uniqueness result for nodal solutions of the standard Ginzburg-Landau model shown in \cite{Hervex2}.
\end{remark}

\subsection{Physical relevance and fine qualitative properties}
\label{intro:physical}

Our analysis of the boundary value problem  \eqref{RS::ODE} \& \eqref{BC} is motivated by  the  study of the energetic stability of the radially-symmetric solution for a system of partial differential equations used for modelling nematic liquid crystals.  This article is the first one in a series of two papers addressing this issue. In the current paper we prove the existence, uniqueness and fine qualitative properties of the radially symmetric solution that is completely determined by the scalar solution $u$ of \eqref{RS::ODE} \& \eqref{BC}, as explained in the remainder of this subsection.  These properties will play an important role in our  second paper \cite{INSZ2} that focuses on proving the energetic stability.

Let us consider  the following energy functional 
 \begin{align}\label{extended energy}
\mcF[Q; \Omega]= \int_{\Omega} \Big[L_1|\nabla{Q}|^2+L_2\nabla_jQ_{ik}\nabla_kQ_{ij}
+L_3\nabla_jQ_{ij}\nabla_kQ_{ik}
+ f_{bulk}(Q)\Big]\,dx,
\end{align} 
where $Q \in H^1(\Omega, \mcS_0), \Omega \subset \RR^3$  with  $$\mcS_0\defeq\{ Q\in \RR^{3\times 3},\, Q=Q^t,\,\tr (Q)=0\}$$ denoting the set of the so-called {\it $Q$-tensors} (here and in the following we assume summation over the repeated indices $i,j,k=1,2,3$).
It is known that the gradient part of the
energy  is bounded from below (and coercive) if
and only if certain relations are assumed between $L_1,L_2, L_3$ (see \cite{Tim-Gart, Longa-several}).
The Euler-Lagrange equations associated to the above energy are:
\begin{multline} \label{def:pde}
2L_1\Delta{Q_{ij}}+(L_2+L_3)\left(\nabla_j\nabla_kQ_{ik}+\nabla_i\nabla_kQ_{jk}\right)-\frac{2}{3}(L_2+L_3)\nabla_l\nabla_kQ_{lk}\delta_{ij}
\\
	=-\left(\frac{\partial f_{bulk}(Q)}{\partial Q}\right)_{ij}+\frac{\delta_{ij}}{3}\tr\left(\frac{\partial f_{bulk}(Q)}{\partial Q}\right), \quad i,j,l,k=1,2,3.
\end{multline} 
In general the bulk potential $f_{bulk}(Q)$ is required to satisfy the physical invariance $f_{bulk}(Q)=f_{bulk}(\Rot Q \Rot^t)$ with $\Rot \in SO(3)$, hence it is a function of the principal invariants of $Q$ (see \cite{ball-isotropic}), which are $\tr(Q^2)$ and $\tr(Q^3)$ (taking into account that $\tr(Q)=0$ in our case). A typical form of the potential often used in the literature is:
\begin{equation}
f_{bulk}(Q) = -\frac{a^2}{2}|Q|^2 - \frac{b^2}{3}\tr(Q^3) + \frac{c^2}{4}|Q|^4,
\label{BulkEDensity}
\end{equation}
where $a^2,c^2>0$, $b^2\ge 0$ and $|Q|^2\defeq\tr(Q^2)$ (see e.g. \cite{Ma-Za} and references therein).

We are interested in studying a radially symmetric solution on balls $\Omega=B_R(0)\subset\RR^3$ with $R\in (0, +\infty]$ (with the convention that $\Omega=\RR^3$ if $R=+\infty$). This solution is relevant in the study of  topological defects in liquid crystals (see \cite{Ma-Za}). More precisely we  say that a matrix-valued measurable map $Q:\Omega\to\mcS_0$  is {\it radially symmetric} if
\begin{equation}
Q(\Rot x) = \Rot \,Q(x)\,\Rot^t \text{ for any } \Rot\in SO(3) \text{ and a.e. } x\in\Omega .
\label{RadSymDef}
\end{equation}
It will be shown in Appendix~\ref{sect:ode_pde} that such a solution of \eqref{def:pde}, called ``the'' {\it melting hedgehog}, can be written as :
\begin{equation}\label{def:MeltingHedgehog}
H(x)=u(|x|)\left(\frac{x}{|x|}\otimes \frac{x}{|x|} -\frac{1}{3}Id\right).
\end{equation} 
In the case of the potential \eqref{BulkEDensity}, $u:\RR_+\to \RR$ -- {\it the scalar profile} of the melting hedgehog -- is a solution of \eqref{RS::ODE} with  $p=2,q=6$ and
\begin{equation}\label{def:physical nonlinearity+}
F(u(r)) =\frac{1}{\alpha}\left(-a^2\,u(r) - \frac{b^2}{3}\,u(r)^2 + \frac{2c^2}{3}\,u(r)^3\right), \quad r>0, 
\end{equation} where $\alpha=2L_1+\frac{4(L_2+L_3)}{3}$, see \cite{Tim-Gart} and \cite{Longa-several}.

As a direct consequence of Theorem~\ref{thm:uniqueness+}, we obtain the following result which is new for the liquid crystal community: the uniqueness of radially symmetric solution of  \eqref{def:pde} is proved in the general class of {\it nodal} scalar profiles.
\begin{theorem}\label{thm:uniquemeltinghedgehog}
Assume that $\alpha=2L_1+\frac{4(L_2+L_3)}{3} > 0$. Consider the equation \eqref{def:pde} with the bulk potential \eqref{BulkEDensity} on the domain $\Omega=B_R(0)$ with the boundary condition~\footnote{The boundary condition is  $\lim_{|x|\to+\infty} \big|Q(x)-s_+\left(\frac{x}{|x|}\otimes \frac{x}{|x|}-\frac{1}{3}Id\right)\big|=0$ if $R=+\infty$ (i.e. $\Omega=\RR^3$).}
$$Q(x)=s_+\left(\frac{x}{|x|}\otimes \frac{x}{|x|}-\frac{1}{3}Id\right)\textrm{ for }x\in\partial B_R(0).$$
Then there exists a unique radially-symmetric solution of the above problem.
 \end{theorem}

One of the important physical questions is related to the stability of this radially symmetric solution as a critical point of the energy \eqref{extended energy}.  Corollary~\ref{cor:energyuniqueness}  shows that the melting hedgehog is locally energy minimizing within the class of radially symmetric tensors, under suitable assumptions on the nonlinearity. The corresponding question of local energy minimality for the melting-hedgehog solution \eqref{def:MeltingHedgehog} with respect to arbitrary perturbations  (with respect to the general energy \eqref{extended energy})  is a considerably more challenging  task and the main motivation for the current work. For the case of physically relevant potential \eqref{BulkEDensity} and $\Omega=\RR^3$, it was shown in \cite{mg} that for $a^2$ large enough the melting hedgehog is not locally stable (hence not locally minimizing) and conjectured, based on numerical evidence, that for $a^2$ small the melting hedgehog is locally stable.  In our forthcoming paper \cite{INSZ2} we prove this conjecture. The crucial step for obtaining the result in \cite{INSZ2} has been a thorough understanding of the  fine qualitative properties of the unique solution $u$ of \eqref{RS::ODE} \& \eqref{BC}. In particular, in \cite{INSZ2} we extensively use  the following result that we prove in Section $4$:

\begin{theorem}\label{thm:ODEqual}
Let $u$ be the unique  solution of \eqref{RS::ODE} and \eqref{BC} where $R=+\infty$, $p=2$, $q=6$ and the right-hand side $F(u)$ is given by \eqref{def:physical nonlinearity}. If we denote $w(r):=\frac{ru'(r)}{u(r)}$ then
\begin{equation}
\label{rel:u'u}
0<w(r) < 2\textrm{ for all }r \in (0,+\infty).
\end{equation} Moreover, setting $f(u)=\frac{F(u)}{u}$, then the following inequalities hold for every $r\in (0, +\infty)$:
\begin{align}
 \label{rel:uprimes}
u'' + \left(-\frac{3u'}{u} + \frac{5}{r}\right)u' &\geq 0,\\
 \label{rel:relfs}
2a^2 + \frac{b^2}{3} u &> -\frac{2}{w}\,f(u),\\
\frac{3}{r^2}(w-2)(w + 1) < f(u)
&< \frac{1}{r^2}(w-2)(2w + 3)  < 0.\label{rel:fws}
\end{align}
\end{theorem}

\subsection{Related literature and organization of the paper}

Let us  now review the existing mathematical literature where similar problems were considered.  The differential equation \eqref{RS::ODE} is a generalization of the equation that describes scalar profiles for Ginzburg-Landau type of equations, as analyzed for instance in \cite{tang, Hervex2, Mironescu-radial}. This type of equations was extensively studied in the last twenty years. Below we mention only few of the papers that are most relevant to our study. 

One of the first results about existence and uniqueness of the solution of Ginzburg-Landau type profile was obtained in \cite{Hervex2}. The authors considered the $2$D case of the Ginzburg-Landau type equation  \eqref{RS::ODE} with the  nonlinearity $F(u)=-u(1-u^2)$ and $p=1$, $q=n^2$ for integers $n\geq 1$. Using shooting method and maximum principle methods they obtained existence and uniqueness of the solution for the problem. The generalization to  higher-dimensional cases was studied in \cite{farina-guedda}, taking $p=n-1$, $q=k(k+n-2)$ for integers $n\ge 3$ and $k\geq 1$. Both papers \cite{Hervex2} and \cite{farina-guedda} investigate {\it nodal} solutions.

For general nonlinearity $F(u)$, existence and uniqueness of {\it positive} solutions are shown in a recent work (see \cite{aguareles-baldoma}) only for the case $p=1$. The authors turn the differential equation into a suitable fixed point equation, and use fixed point methods and a sliding method to show existence and uniqueness of the positive solution. Moreover they also obtain some results on a qualitative behavior of the solution. 

The profile of the radially symmetric solution for Landau-de Gennes problem has been recently studied in \cite{lamy}. Using Pohozaev-type arguments  the author showed the monotonicity and uniqueness of the energy-minimizing solution of equation  \eqref{RS::ODE} in bounded domains for $F(u)$ of type \eqref{def:physical nonlinearity}.

In this paper we consider the equation \eqref{RS::ODE} with $p, q \in \RR$, $q>0$ and general nonlinearity $F(u)$ on bounded and unbounded domains. We show existence and uniqueness of positive solutions with very light and natural restrictions on $F(u)$. 
Moreover, we also show uniqueness of general nodal solutions for $p\geq 0$ under more restricted assumptions on nonlinearity $F(u)$.  Using the mountain pass theorem, we provide a counterexample to uniqueness of nodal solution when $F(u)$ does not satisfy these assumptions.
Finally, we investigate fine properties of the solution corresponding to the radially symmetric profile of the melting hedgehog in Landau-de Gennes model of liquid crystals. These fine properties are of utter importance in the investigation of the stability of the melting hedgehog  that we perform in the forthcoming paper \cite{INSZ2}.
  
  \bigskip
The paper is organized as follows.  In Sections \ref{sec:exist} and \ref{sec:uniqueness} we gather the arguments for proving Theorem~\ref{thm:main} on the existence and uniqueness of positive solutions to \eqref{RS::ODE} \& \eqref{BC}. The proof of Theorem~\ref{thm:main} is provided at the end of Section~\ref{sec:uniqueness}. Corollary~\ref{cor:energyuniqueness} on locally energy minimizing solutions is shown in Section~\ref{ssec:NoPosAs}, where we also prove Theorem~\ref{thm:uniqueness+} on the uniqueness of nodal solutions. Theorem~\ref{thm:ODEqual} is proved in  Section \ref{sec:refined} where certain refined properties of the solution corresponding to the nonlinearity \eqref{def:physical nonlinearity} are studied. Section \ref{sec:SC} is devoted to proving the existence of a sign-changing solution of \eqref{RS::ODE} \& \eqref{BC} for certain types of nonlinearities (see Proposition~\ref{lem:MP}). In Appendix \ref{sect:ode_pde} we provide some properties of radially symmetric $Q$-tensors. Finally, in Appendix \ref{app} we present versions of maximum principle that are needed in the body of the paper.


\section{Existence and behaviour near $0$ and $+\infty$}
\label{sec:exist}

In this section we prove the existence of solutions of the problem \eqref{RS::ODE}\&\eqref{BC} under \eqref{cond:pq}. When $R$ is finite, this is done via an energy minimization procedure. The case $R = +\infty$ is obtained by a limiting process. A delicate issue will be to ensure that the solution thus obtained in the limit does not become trivial and has the desired asymptotic behaviours at $0$ and $+\infty$.

\subsection{Existence on finite domains}
\label{ssec:ExistenceFinite}

For $F:\RR\to\RR$ with \eqref{condF} we associate $\tilde F:\RR\to\RR$ to  be any $C^1$ function such that
\begin{equation}\label{FCond:even}
\tilde F(t) = F(t) \text{ for }t\ge 0\textrm{ and } \tilde F_{even}(t)\textrm{ satisfies }\eqref{Fcond:even}.
\end{equation}
(For example, we can define $\tilde F$ by $\tilde F(-t) = - F(t)$ for $t > 0$.) Let
\[
\tilde h(t) = 2\int_0^t \tilde F(s)\,ds.
\]
Note that, by \eqref{FCond:even}, we have
\begin{equation}\label{hAsymm}
\tilde h(-|t|) \ge \tilde h(|t|)  \text{ for all } t \in \RR,
\end{equation}
and so, by \eqref{condF}, $\tilde h$ is bounded from below.

Consider instead of the energy $E$ defined by \eqref{energyODE} the following modified energy:
\begin{equation}\label{energyODE0}
\tilde E[u;(0,R)]= \frac{1}{2}\int_0^R \Big[r^{p}|u'(r)|^2+qr^{p-2}u(r)^2+r^p \tilde h(u(r))\Big]\,dr.
\end{equation}
Since $\tilde F \equiv F$ in $[0,+\infty)$, all non-negative critical points of $\tilde E$ coincide with non-negative critical points of $E$ and vice versa, as can be seen by looking at the corresponding Euler-Lagrange equations. In other words, if we are interested in positive solutions of \eqref{RS::ODE} we can always assume that $F$ satisfies \eqref{Fcond:even}.

\begin{lemma}\label{lemma:uR} 
Assume \eqref{cond:pq}, \eqref{condF} and \eqref{FCond:even}. Then for every $R \in (0,+\infty)$, there exists a global energy minimizer $u_R$ of $\tilde E$ over 
\[
\mcM_R:=\Big\{u:(0,R)\to\RR\, :\, r^{p/2}u', r^{p/2-1}u\in L^2(0,R),\quad u(R)=s_+\Big\}.
\]
Moreover $u_R$ satisfies \eqref{RS::ODE} and $0\le u_R(r)\le s_+$ for all $r\in (0,R)$.
\end{lemma}

\begin{proof} We split the proof in several steps.

\bigskip

\medskip\noindent{\it Step 1: Reduction from $\mcM_R$ to $\mcM'_R$ where}
$$\mcM'_R:=\{u\in \mcM_R\, :\, 0\le u(r)\le s_+, \, r\in (0, R)\, \} \subset \mcM_R.$$
We claim that $$\inf_{\mcM_R}\tilde E =  \inf_{\mcM'_R}\tilde E.$$
To this end let us take $u\in \mcM_R\setminus\mcM'_R$. Set $$\bar u(r)= |u|(r), \, r\in (0,R).$$ Then $\bar u\in\mcM_R$, and by \eqref{hAsymm},
\begin{equation}\label{red1}
\tilde E[\bar u;(0,R)] \leq \tilde E[u;(0,R)].
\end{equation}
We  define now $$\tilde u(r)= \min(\bar u(r),s_+).$$ Then $\tilde u \in \mcM'_R$ and thanks to the fact that $q > 0$ and $\tilde h'(t)=2\tilde F(t)=2F(t)\geq 0$ for $t>s_+$ (by \eqref{condF}) so $\tilde h(\bar u)\ge \tilde h(\tilde u)$ in $(0, R)$, we have 
\begin{equation}\label{red2}
\tilde E[\tilde u;(0,R)] \leq \tilde E[\bar u;(0,R)].
\end{equation} 
The claim follows from \eqref{red1} and \eqref{red2}.

\bigskip

\medskip\noindent{\it Step 2: $\inf_{\mcM'_R} \tilde E >-\infty$.} Indeed, if $p > -1$, this step is clear since $\tilde h=h$ is bounded in the interval $[0,s_+]$ and the function $r\mapsto r^p$ is integrable on $(0,R)$. In the general case, for $p\in\RR$, we argue as follows. Since $F(0)=0$ and $|F'|\leq C_1$ on $[0, s_+]$ with $C_1>0$, we have $|F(t)|\le C_1t$ for  $0 \leq t \le s_+$. Hence 
\begin{equation}\label{est:hquadratic}
|h(t)|\le C_1|t|^2\,\textrm{ for } t \in [0,s_+].
\end{equation} 
Moreover, by \eqref{condF}, we have for  $t\in [0,s_+]$ :
\begin{equation}
\label{est:hLow}
0\ge h(t)=2\int_0^t F(s)\,ds \ge 2\int_0^{s_+} F(s)\,ds.
\end{equation}
Set $u\in {\mcM'_R}$. For $0 < r \leq R_0 = (\frac{q}{2C_1})^{1/2}$, by \eqref{est:hquadratic}, we have $r^ph(u(r)) \ge -\frac{q}{2}r^{p-2}u(r)^2$, while for $R_0 \leq r \leq R$, we have by \eqref{est:hLow}:  $$r^ph(u(r)) \ge 2\max(R_0^p,R^p)\int_0^{s_+} F(s)\,ds=:-C_2,$$ 
with $C_2>0$. It follows that
\begin{equation*}
r^ph(u(r))\ge -\frac{q}{2} r^{p-2}u^2(r)-C_2,\forall r\in (0,R).
\end{equation*}
Thus, the function
\begin{equation}\label{lowerbound:rphu}
T(u)(r):=r^p h(u(r))+qr^{p-2} u^2(r)+C_2\ge 0,\, r\in (0, R)\end{equation}
is positive and therefore, we have $\tilde E(u)\ge - C_2 R/2>-\infty$ for every $u\in {\mcM'_R}$, which finishes Step 2. Note that $\inf_{\mcM'_R} \tilde E <\infty$ since every configuration $u\in {\mcM'_R}$ with $u\equiv 0$ near $r=0$ has finite energy $\tilde E(u) <\infty$.

\bigskip

\medskip\noindent{\it Step 3: Existence of a minimizer of $\tilde E$ over $\mcM'_R$.}
Indeed, by \eqref{lowerbound:rphu},
the direct method of calculus of variation using Sobolev's embedding and Fatou's lemma establishes the existence of a minimizer of $\tilde E$ over $\mcM'_R$. We omit the details.
\end{proof}

\begin{remark}
Let us point out that since the potential $\tilde F$ satisfies the condition \eqref{Fcond:even}, we can use the uniqueness result given by Corollary~\ref{cor:energyuniqueness} (to be proved in the next section) and show that
$\textrm{argmin}_{\mcM_R} \tilde E=\textrm{argmin}_{\mcM'_R}\tilde E$ and it contains one single element.
\end{remark}

To complete the proof of the existence in the case of a finite domain, we need to show that $u_R(0) = 0$. In fact, we prove stronger asymptotic estimates in the next subsection.


\subsection{Local behaviour near the origin}
\label{ssec:LocBeZero}

Note that the homogeneous linear equation associated with \eqref{RS::ODE} is a Fuchsian ODE at $r = 0$, see e.g. \cite{Birkhoff-Rota}. Let $\gamma_\pm$ denote the solutions of the indicial equation, i.e.
\be
\label{Fuch_ind}
\gamma_\pm:=\frac{1-p\pm\sqrt{(p-1)^2+4q}}{2}.\ee
As $q > 0$, we have that $\gamma_+ > 0 > \gamma_-$. Thus, if $u$ is a bounded solution of \eqref{RS::ODE}, then we  expect that $u$ ``behaves like $r^{\gamma_+}$'' at the origin.

\begin{proposition}\label{Prop:OAsymp}
Assume that condition \eqref{cond:pq} holds and  $F$ is a $C^1$ function satisfying $F(0) = 0$. Let $u$ be a (nodal) solution of \eqref{RS::ODE}
on $(0, R)$ with  $R\in (0, +\infty]$ such that $u$ is bounded near the origin.

{\it (i)} Then the function $v(r) := \frac{u(r)}{r^{\gamma_+}}$ is differentiable up to $0$ and $v'(0)= 0$. In particular $u(0)=0$.

{\it (ii)} If in addition, $F$ satisfies  \eqref{condF} and $u\ge 0$ in $(0, R)$ and $u(R)\in (0, s_+]$, then $0<u<s_+$ on $(0,R)$. Moreover $v$ is decreasing and in particular
\begin{equation}\label{bound:derivgeneral} 
u'(r) < \frac{\gamma_+ u(r)}{r} \;\textrm{ on }(0,R).
\end{equation}
\end{proposition}

\bigskip

Note that if $F$ satisfies  \eqref{condF}  and the first condition in \eqref{condFLeft}, then every solution $u$ of \eqref{RS::ODE} with $u(R)=s_+$
is bounded (i.e., $s_-\leq u\leq s_+$ in $(0,R)$ by the maximum principle) and therefore, Proposition \ref{Prop:OAsymp} implies that
$u$ satisfies \eqref{BC}. 

\bigskip

\begin{proof} Assume that $|u(r)| \leq M$ for $r \in (0,\delta_0)$ for some $\delta_0 \in (0, R]$. Standard regularity result for ODEs implies that $u\in C^3(0,R)$.

\medskip\noindent{\it Step 1: We first show that 
\begin{equation}
|u| \leq C\,r^{\gamma^+} \text{ in }(0,\delta_0)
\label{Eq:28IV13.01}
\end{equation}
with $C>0$ depending only on $M$.} Indeed, denoting $u_\pm:=\max\{0, \pm u\}$, we prove \eqref{Eq:28IV13.01} for both $u_\pm$. Since $F$ is $C^1$ with $F(0) = 0$, we have $|F(t)| \leq \tilde C|t|$ for $t \in [-M,M]$ with the constant $\tilde C>\|F'\|_{L^\infty(-M,M)}\geq 0$ depending only on $M$; in particular, 
$$F(-u_-)\leq \tilde C u_- \quad \textrm{and} \quad -F(u_+)\leq \tilde C u_+\quad \textrm{ in } (0, \delta_0).$$
Then, by \eqref{RS::ODE}, we deduce:
$$Lu_+:=-u_+''-\frac pr u_+'+\big(\frac{q}{r^2}-\tilde C\big)u_+\le 0 \quad \textrm{ as measure in } (0,  \delta_0).$$ By theory of ODEs with a regular singular point (see for instance \cite{Birkhoff-Rota}), there exist functions $w_1$, $w_2$ such that $
w_1(r) = r^{\gamma_+} + o(r^{\gamma_+})$ and $w_2(r) = r^{\gamma_-} + o(r^{\gamma_-})$ as $r\to 0$, and $Lw_1=Lw_2=0$. By choosing $\delta_0>0$ smaller if necessary we can assume that 
$w_1,w_2> 0$ on $(0,\delta_0]$ and $\frac{q}{r^2}\geq \tilde C $ on $(0,\delta_0]$.
Note now that for any constant $\mu$ it holds:
$$L(\mu\omega_1-u_+)\ge 0=Lw_2.$$
Choosing $\mu>0$  such that $\mu\omega_1(\delta_0)\ge M\ge u_+(\delta_0)$ we can apply Lemma~\ref{lemma:weakmax} and obtain $\mu\omega_1\ge u_+$ on $(0,\delta_0)$. 
The estimate for $u_-$ follows by the same argument (since $Lu_-\leq 0$ in the sense of measures on $(0,  \delta_0)$). 
Noting that $\mu$ depends only on $M$ (and not on $u$), we obtain the claimed \eqref{Eq:28IV13.01}. In particular, $u(0)=0$.

\medskip\noindent{\it Step 2: We prove that $v$ is differentiable up to $r = 0$.}
In view of \eqref{Eq:28IV13.01}, we deduce from \eqref{RS::ODE} that
\[
\Big|u'' + \frac{p}{r} u' - \frac{q}{r^2} u\Big|=|F(u)|\le \tilde C|u|\leq \bar C\,r^{\gamma_+} \text{ in } (0,\delta_0).
\]
Denote $$L_0u=-u''-\frac pr u'+\frac{q}{r^2} u.$$ Then we have
$$-\bar Cr^{\gamma_+}\le L_0u\le \bar Cr^{\gamma_+} \text{ in } (0,\delta_0).$$
Note that $L_0 r^{\gamma_\pm}=0$ and $L_0(r^{\gamma_++2})=-2(2\gamma_++p+1)r^{\gamma_+}$ with $2\gamma_++p+1>0$. Set $\tilde u^\pm:=u\pm \frac{\bar C}{2(2\gamma_++p+1)}r^{\gamma_++2}$. Then $L_0 \tilde u^+\leq 0$ and $L_0 \tilde u^-\geq 0$  in $(0, \delta_0)$. 
Let $s\in (0, \delta_0)$ and note that
we are in the framework of Lemma~\ref{lemma:weakmax}  (with $w_0=r^{\gamma_-}$) applied to 
$$L_0(\mu_\pm r^{\gamma_+}\mp\tilde u^\pm)\ge 0 \, \textrm{ on }(0, s)$$ where $\mu_\pm=\mu_\pm(s)\in \RR$ is determined by $\mu_\pm s^{\gamma_+}:=\pm \tilde u^\pm(s)
$. We deduce:
$$\pm \tilde u^{\pm} (r)\le \frac{\pm \tilde u^{\pm} (s)}{s^{\gamma_+}}r^{\gamma_+}, \quad  0<r<s.$$ It follows that
\be\label{est:vgrad}
\big|\frac{u(r)}{r^{\gamma_+}}-\frac{u(s)}{s^{\gamma_+}}\big|\le O(s^2-r^2)
\ee for  $0<r<s$. Since $s$ was arbitrarily chosen in $(0, \delta_0)$, we have that \eqref{est:vgrad} implies the existence of a limit of $v$ at the origin. Dividing \eqref{est:vgrad} by $s-r$ and passing to the limit $r\to 0$, followed by $s\to 0$, we obtain $v'(0)=0$. Since $u$ is $C^2$ away from $0$, we conclude that $v$ is differentiable up to the origin which ends the proof of ${\it (i)}$. 

\medskip
\noindent{\it Step 3: Proof of (ii)}. Assume that the stronger hypothesis in $(ii)$ holds. First, by \eqref{condF}, we note that $L_0 u\geq 0$ in $(0, R)$ and $u\geq 0$ in $(0, R)$; thus, by the strong maximum principle, if $u$ achieves the value $0$ inside the interval $(0,R)$, it must be identically  zero which would violate $u(R)>0$. So, $u>0$ in $(0, R)$.
Second, note that $u \leq s_+$ in $(0,R)$ because otherwise, $u$ would achieve a local maximum at some $r_0\in (0, R)$ where $u(r_0)>s_+$ and
\[
0 >-L_0 u(r_0)= u''(r_0) + \frac{p}{r_0}u'(r_0) - \frac{q}{r_0^2}\,u(r_0) = F(u(r_0)) \stackrel{\eqref{condF}}{\geq}0,
\]
which is absurd. Third, note that $L_0 u+F(u)=0\leq L_0 (s_+)+F(s_+)$. Therefore, we obtain $$M(s_+-u):=L_0 (s_+-u)+a(r)(s_+-u)\geq 0$$ where
$M$ is a linear elliptic operator with
 $a$ a bounded continuous function defined by $a(r)=\frac{F(s_+)-F(u(r))}{s_+-u(r)}$ if $u(r)\neq s_+$ and $a(r)=F'(s_+)$ otherwise. As above, the strong maximal principle applied for $M$ and $s_+-u\geq 0$ implies that $u<s_+$ on $(0,R)$ (because of $u(0)=0$ which prevents $u$ being identically constant to $s_+$).

It remains to show that $v$ decreases.
For that, note first that $v$ satisfies 
\begin{equation}
(r^{2\gamma_++p} v')' = r^{2\gamma_+ + p}\Big(v'' + \frac{2\gamma_+ + p}{r} v'\Big) = r^{\gamma_+ + p}\Big(u'' + \frac{p}{r} u' - \frac{q}{r^2} u\Big) = r^{\gamma_+ + p} F(u).
\label{15V13-E1}
\end{equation}
Using \eqref{condF}, we obtain $(r^{2\gamma_++p} v')' < 0$ on $(0,R)$ because $0<u<s_+$ on $(0,R)$, meaning that $r^{2\gamma_++p} v'$ is decreasing on $(0,R)$.
Noting that $2\gamma_++p>0$ and $v'(0)=0$ (by step 2), it follows:
\begin{equation}
\lim_{r \rightarrow 0} r^{2\gamma_+ + p}v'(r) =   0.
	\label{15V13-E2}
\end{equation}
Therefore, it follows $r^{2\gamma_++p} v'< 0$ on $(0,R)$. So we conclude that $v' < 0$, i.e. $v$ is decreasing on $(0,R)$. Estimate \eqref{bound:derivgeneral} is now straightforward.
\end{proof}

\begin{corollary}
\label{coro_finit}
Assume \eqref{cond:pq} and \eqref{condF}. For every $R\in (0,+\infty)$, there exists a non-negative solution $u$ of \eqref{RS::ODE}\&\eqref{BC} that is a minimizer of the energy $E$ defined in \eqref{energyODE} over the set of non-negative configurations $\{v\in \mathcal{M}_R\, :\, v(r)\geq 0, \, r\in (0,R)\}$ where $\mathcal{M}_R$ is defined in Lemma \ref{lemma:uR}. Moreover, $0< u < s_+$ in $(0,R)$
and $u$ is increasing on $(0,R)$.
\end{corollary}

\begin{proof} The first part of the statement is a direct consequence of Lemma \ref{lemma:uR} and Proposition \ref{Prop:OAsymp} since the energy $E$ coincides with the energy $\tilde E$ for non-negative configurations in $\mathcal{M}_R$. The fact that $u$ is increasing is a consequence of Lemma \ref{ERS:monotonicity} that we postpone the proof for Section \ref{sec:uniqueness}.
\end{proof}


\subsection{Existence on infinite domain}

Let us now prove the existence of solution to \eqref{RS::ODE}\&\eqref{BC} in case $R=+\infty$. 

\begin{proposition}\label{prop:existenceinftydomain}
Assume \eqref{cond:pq} and \eqref{condF}. For $R=+\infty$, there exists a non-negative increasing solution $u$ to \eqref{RS::ODE}\&\eqref{BC}. Furthermore, $0 < u < s_+$ in $(0,R)$ and $u$ is locally minimizing with respect to the energy $\tilde E$ defined in \eqref{energyODE0}.
\end{proposition}

\begin{proof}  We  proceed in several steps:

\medskip\noindent{\it Step 1: Constructing a solution of \eqref{RS::ODE} on $(0,+\infty)$} We denote by $u_n$ a global energy minimizer of the energy $\tilde E$ obtained in Lemma~\ref{lemma:uR} on the interval $(0,n)$ and in the space $\mcM_n$ that satisfies $u_n \in [0,s_+]$. We extend $u_n$ to the function $\bar u_n$ on $[0,+\infty)$ by letting $\bar u_n(r)=\left\{\begin{array}{ll} u_n(r) &\textrm{ if }r\in (0,n)\\ s_+ &\textrm{ if }r>n\end{array}\right.$. 
Obviously, the sequence $(\bar u_n)_{n\in\NN}$ is uniformly bounded in $L^\infty(0,+\infty)$. Let $I\subset (0,+\infty)$ be a compact interval and $n_0\in \NN$ so that
$I\subset (0, n_0)$. By standard regularity arguments for the ODE \eqref{RS::ODE}, one can show that $(\bar u_n)_{n\geq n_0}$ is uniformly bounded on $C^3(I)$. Since $I$ is arbitrarily chosen, by Arzela-Ascoli's theorem, we deduce that $(\bar u_{n})$ converges (up to a subsequence) in $C^2_{loc}(0,+\infty)$ to some $u_\infty\in C^2(0,\infty)$ which satisfies \eqref{RS::ODE} and $u_\infty \in [0,s_+]$.

\medskip\noindent{\it Step 2: Behaviour of $u_\infty$ at $0$.} Since $u_\infty$ satisfies \eqref{RS::ODE} and $u_\infty \in [0,s_+]$, Proposition \ref{Prop:OAsymp} implies that $\frac{u_\infty}{r^{\gamma_+}}$ is differentiable up to the origin. In particular $u_\infty(0) = 0$.

\medskip\noindent{\it Step 3: Behaviour of $u_\infty(r)$ as $r\to+\infty$.} We know that $\bar u_n$ are non-decreasing functions on $(0,+\infty)$ by Corollary \ref{coro_finit}. Then the limit function $u_\infty$ is also non-decreasing. Since $0\le u_\infty\le s_+$, then there exists 
$$s_\infty:=\lim_{r\to+\infty}u_\infty(r) \in [0,s_+].$$
{\it Claim: $s_\infty \in \{0,s_+\}$.} 

\begin{proof} Assume by contradiction that $0<s_\infty<s_+$. Recall: 
$$\frac{1}{r^p}(r^p u_\infty')'=q\frac{u_\infty}{r^2}+F(u_\infty).$$ As $r\to+\infty$ we have $u_\infty(r)\to s_\infty, F(u_\infty(r))\to F(s_\infty)<0$; hence, for $\varepsilon>0$ small enough there exists $R_0>0$ so that 
$$\frac{1}{r^p}(r^pu_\infty')'\le -\varepsilon \quad \textrm{for $r\ge R_0$}.$$ 
If $p=-1$, we integrate the above inequality on $(R_0, r)$ to obtain:
$$\frac{u_\infty'(r)}{r}\le \frac{u_\infty'(R_0)}{R_0}-\varepsilon(\log r-\log R_0)\to -\infty \quad \textrm{as $r\to+\infty$}.$$ We deduce that $u'_\infty(r)<0$ for $r$ large enough, 
which contradicts the fact that $u_\infty$ is non-decreasing. Consider now $p\not=-1$. As before, integrating on $(R_0, r)$, we obtain:
\begin{equation}\label{rel:rpsinfty}
r^p u_\infty'(r)\le R_0^p u_\infty'(R_0)-\frac{\varepsilon}{p+1}(r^{p+1}-R_0^{p+1}).
\end{equation} 
We have now two cases:

\noindent{\it Case $p>-1$.} As before, $r^{p+1}\to+\infty$ as $r\to+\infty$ and \eqref{rel:rpsinfty} implies $r^p u_\infty'(r)<0$ for $r$ large enough, obtaining again a contradiction.

\smallskip\noindent{\it Case $p<-1$.} Relation \eqref{rel:rpsinfty} implies

$$u_\infty'(r)\le \frac{R_0^{p+1}}{r^p}\left(\frac{u_\infty'(R_0)}{R_0}+\frac{\varepsilon}{p+1}\right)-\frac{\varepsilon}{p+1}r.$$ By Proposition \ref{Prop:OAsymp}, we deduce that $u_\infty'(r) < \frac{\gamma_+ u_\infty(r)}{r}
\leq \frac{\gamma_+ s_+}{r}$ on $(0,R)$.
Therefore, we choose now $R_0$ large enough such that
$$0\le \frac{u_\infty'(R_0)}{R_0}\le -\frac{\varepsilon}{2(p+1)}.$$
Then $\frac{R_0^{p+1} }{r^p}\left[ \frac{\varepsilon}{2(p+1)}\right] - \frac{\varepsilon}{p+1}r\to -\infty$ as $r\to+\infty$ and we obtain again $u'_\infty(r)<0$ for $r$ large enough, which contradicts the fact that $u_\infty$ is non-decreasing. In all the cases, we obtain that $s_\infty\in \{0,s_+\}$ which concludes the Claim.
\end{proof}

\medskip\noindent{\it Step 4: $u_\infty$ is locally minimizing w.r.t. energy $\tilde E$.} 
Let $\omega\subset (0,+\infty)$ be a compact interval and $n_0\in \NN$ so that $\omega\subset (0, n_0)$. 
Since $u_n$ is a global minimizer for $\tilde E[\cdot;(0,n)]$, we have for any $n\geq n_0$ that $\tilde E[u_n;\omega]\le \tilde E[u_n+\varphi;\omega]$ for any $\varphi\in C_c^\infty(\omega)$. As $u_n\to u_\infty$ in $C^2(\omega)$, we can pass to the limit in the above inequality and obtain that $u_\infty$ is locally energy minimizing.

\medskip\noindent{\it Step 5: Showing that $u_\infty\not\equiv 0$ and $s_\infty = s_+$.} We assume by contradiction that $u_\infty\equiv 0$. Since it is locally minimizing, we have 
for any compact interval $\omega\subset (0,+\infty)$ that
\begin{equation}\label{0min}
E[0;\omega]=0\le \tilde E[\varphi;\omega] \text{ for any } \varphi\in C_c^\infty(\omega). 
\end{equation} 
Let us pick an arbitrary $\varphi\in C_c^\infty (0,1)$ with $\varphi\not\equiv 0$ and $\varphi\in [0,s_+]$. Set $\varphi_n(r):=\varphi(\frac{r}{n})$ for every $r>0$ so that 
$\f_n\in C_c^\infty (0,n)$.
We have
\begin{align*}
2\tilde E[\varphi_n;(0,n)]=2E[\varphi_n;(0,n)]&=\int_0^n \left[\frac{r^p}{n^2}(\varphi'(\frac{r}{n}))^2+qr^{p-2}\varphi^2(\frac{r}{n})+r^p h(\varphi(\frac{r}{n}))\right]\,dr\\
&= n^{p+1}\left[ \frac{1}{n^2}\bigg(\int_0^1 \big(t^p (\varphi')^2+qt^{p-2}\varphi^2\big)\,dt\bigg)+\int_0^1 t^p h(\varphi)\,dt\right].
\end{align*}
However, $\int_0^1 h(\varphi(t))\,dt<0$ as $h<0$ on $(0,s_+)$, so
$$\frac{1}{n^2}\bigg(\int_0^1 \big(t^p (\varphi')^2+qt^{p-2}\varphi^2\big)\,dt\bigg)<\left|\int_0^1 t^p h(\varphi)\,dt\right| \text{ for $n$ large enough}.$$ This implies $E[\varphi_n;(0,n)]<0$ for n large which contradicts \eqref{0min}. So, $u_\infty\not\equiv 0$. Since $u_\infty$ is non-decreasing, it means that $s_\infty>0$ so that by Step 3, we conclude that $s_\infty=s_+$.

Finally, by Proposition \ref{Prop:OAsymp}, we deduce that $u_\infty\in (0, s_+)$ on $(0, R)$ and by Lemma \ref{ERS:monotonicity}, we conclude that $u_\infty$ is increasing on $(0, R)$.
\end{proof}

\subsection{Local behaviour near infinity}
\label{ssec:LocBeInfty}

On infinite domains, we study the asymptotic behavior of a solution $u$ near $R=+\infty$.

\begin{proposition}\label{ERS::Asymp}
Assume \eqref{cond:pq} and \eqref{condF}. If $u$ is a non-negative solution of (\ref{RS::ODE}) \&(\ref{BC}) with $R = +\infty$, then 
\begin{equation}\label{asympt:u}
u(r) = s_+ -\frac{\beta}{r^2} + o(r^{-2}) \text{ as } r\rightarrow +\infty
	\;,
\end{equation}
where
\begin{equation}\label{eq:beta}
\beta=\frac{q \,s_+}{F'(s_+)}
	.
\end{equation}

\end{proposition}
\begin{remark}
If we assume (\ref{asympt:u}) then the value of $\beta$ in \eqref{eq:beta} can be formally computed by matching the powers of $\frac{1}{r^2}$ in (\ref{RS::ODE}) as $r\to +\infty$. 
\end{remark}

\begin{proof} We divide the proof in several steps.

\nd {\it Step 1. A change of variables.} Define $$\psi(\tau)=s_+-u(\frac{1}{\tau}),\quad \forall \tau>0.$$ Then $\psi(\tau)\in (0, s_+)$ for $\tau>0$ since $u\in (0, s_+)$ on $(0, R)$ by Proposition \ref{Prop:OAsymp}. A straightforward computation shows that $\psi(\tau)$ satisfies the equation
\begin{equation}\label{eq:psi}
-\psi''(\tau)+\frac{p-2}{\tau}\psi'(\tau)+\frac{q}{\tau^2}\psi(\tau)=\frac{qs_+}{\tau^2}-\frac{F'(s_+)+z(\psi(\tau))}{\tau^4}\psi(\tau), \quad \tau>0,
\end{equation} 
where $$z(s):=\frac{F(s_+)-F(s_+-s)}{s}-F'(s_+), \quad s\in (0, s_+).$$ Obviously, $\lim_{s\to 0} z(s)=0$ and $\lim_{\tau\to 0} \psi(\tau)=0$ (by \eqref{BC}). We will prove that
$\psi(\tau)/\tau^2$ converges as $\tau \to 0$.
\medskip

\nd {\it Step 2. Upper bound of $\psi(\tau)/\tau^2$.} We denote
\begin{equation}\label{def:epsilondelta}
\varepsilon(\delta):=\max_{\tau\in [0,\delta]} \big|q\tau^2+z(\psi(\tau))\big|.
\end{equation}
By Step 1, we have $\lim_{\delta\to 0}\varepsilon(\delta)=0$. Then, by \eqref{condF}, there exists a $\delta_0>0$ so that 
\begin{equation}\label{fixing delta}
|6-2p|<\frac{F'(s_+)-\varepsilon(\delta)}{\delta^2},\quad \forall  \delta\in (0,\delta_0).
\end{equation} 
Fix now $\delta\in (0,\delta_0)$ and set
$$L\psi:=-\psi''+\frac{p-2}{\tau}\psi'+ \frac{\psi}{\tau^4}(F'(s_+)-\varepsilon(\delta)).$$
Then (\ref{eq:psi}),(\ref{def:epsilondelta}) and $\psi\ge 0$ imply that
\begin{equation}
L\psi(\tau)\le\frac{qs_+}{\tau^2},\forall 0<\tau<\delta.
\label{relat_lpsi}
\end{equation}
An upper bound on $\psi$ is provided by means of a suitable comparison function and the weak maximum principle  in Lemma~\ref{lemma:weakmax} applied to $L$ (see also Remark~\ref{remark:fundsol}). We take  
\be
\label{defD}
\phi(\tau)=D\tau^2\quad \textrm{ for } \quad D = D(\delta) :=\max\{\frac{qs_+}{F'(s_+)-\varepsilon(\delta)-|6-2p|\delta^2},\frac{\psi(\delta)}{\delta^2}\}>0.\ee 
Then, by (\ref{fixing delta}),
$$L\phi(\tau)=  D(-6+2p+\frac{F'(s_+)-\varepsilon(\delta)}{\tau^2})\ge \frac{qs_+}{\tau^2}\stackrel{\eqref{relat_lpsi}}{\geq} L\psi(\tau), \quad \forall 0<\tau<\delta.$$
Also, by (\ref{defD}), $\phi(\delta)\ge \psi(\delta)$. The weak maximum principle in Lemma \ref{lemma:weakmax} applied to the operator $L$ and $(\phi-\psi)$ implies that
\[
\psi(\tau)\le \phi(\tau) = D\tau^2, \quad \forall \tau\in (0,\delta).
\]

\medskip

\nd {\it Step 3. Lower bound of $\psi(\tau)/\tau^2$.}  Analogously, we have
\[
\tilde L\psi(\tau):=-\psi''+\frac{p-2}{\tau}\psi'+ \frac{\psi}{\tau^4}(F'(s_+)+\varepsilon(\delta)) \ge\frac{qs_+}{\tau^2}, \quad \forall \tau\in (0,\delta).
\]
Thus, if we denote $\tilde\phi(\tau)=\tilde D\tau^2$ with $\tilde D = \tilde D(\delta) :=\min\{\frac{qs_+}{F'(s_+)+\varepsilon(\delta)+|6-2p|\delta^2},\frac{\psi(\delta)}{\delta^2}\}$ so that
\[
\tilde L\tilde\phi(\tau)=-2\tilde D+2\tilde D(p-2)+\frac{\tilde D(F'(s_+)+\varepsilon(\delta))}{\tau^2}<\frac{qs_+}{\tau^2}, \quad \forall \tau\in (0,\delta),
\]
then we can apply Lemma \ref{lemma:weakmax} to arrive at
\[
\psi(\tau)\ge \tilde\phi(\tau) = \tilde D\tau^2,\quad \forall \tau\in (0,\delta).
\]
Together with Step 2, we conclude that
\begin{equation}\label{ineq:psicomparison}
\tilde D(\delta)\le\frac{\psi(\tau)}{\tau^2}\le D(\delta), \quad \textrm{ for all }  0 < \tau < \delta < \delta_0.
\end{equation}

\medskip

\nd {\it Step 4. We prove that the limit $\lim_{\tau\to 0}\frac{\psi(\tau)}{\tau^2}$ exists.} We denote 
$$\underline{\beta}:=\liminf_{\tau\to 0} \frac{\psi(\tau)}{\tau^2}\textrm{ and }\overline{\beta}:=\limsup_{\tau\to 0}\frac{\psi(\tau)}{\tau^2}.$$ 
We let $p_k\to 0$ and $P_k\to 0$ be sequences such that $\lim_{k\to 0}\frac{\psi(p_k)}{p_k^2}=\underline{\beta}$ and $\lim_{k\to 0}\frac{\psi(P_k)}{P_k^2}=\overline{\beta}$. We  assume without loss of generality that $p_{k+1} < P_k<p_k,\forall k\in\mathbb{N}$. Replacing $\tau=P_k$ and $\delta=p_k$ in (\ref{ineq:psicomparison}) and letting $k\to +\infty$ we obtain:
\begin{equation}\label{ineq:rightPp}
\overline{\beta}\le \max\{\frac{qs_+}{F'(s_+)},\underline{\beta}\}.
\end{equation} 
Likewise we have
\begin{equation}\label{ineq:leftPp}
\min\{\frac{qs_+}{F'(s_+)},\overline{\beta}\}\le \underline{\beta}.
\end{equation}
One can easily see that  (\ref{ineq:rightPp}) and (\ref{ineq:leftPp}) imply $\underline{\beta}=\overline{\beta}$ thus proving our claim that the limit $\beta:=\lim_{\tau\to 0}\frac{\psi(\tau)}{\tau^2}$ exists.

\medskip

\nd {\it Step 5. We prove \eqref{eq:beta}.} If we know in addition that $\tau^2 \psi''$ or $\tau \psi'$ converges to zero as $\tau \rightarrow 0$, \eqref{eq:beta} can be derived immediately from \eqref{eq:psi}. Since we do not assume such convergence, we proceed as follows. Let us denote $\tau_k:=2^{-k}$ and observe that by mean value theorem there exists $\sigma_k\in (\tau_{k+1},\tau_k)$ so that 
\begin{equation}\label{eq:psisigma}
\psi'(\sigma_k)=\frac{\psi(\tau_k)-\psi(\tau_{k+1})}{\tau_k-\tau_{k+1}}=2\tau_k\frac{\psi(\tau_k)}{\tau_k^2}-\tau_{k+1}\frac{\psi(\tau_{k+1})}{\tau_{k+1}^2}   \to 0 \quad \textrm{as } \, k\to +\infty,
\end{equation}
where we used Step 4. We multiply (\ref{eq:psi}) by $\tau^2$, integrate over $[\sigma_{k+2},\sigma_k]$ and by parts, obtaining:
\begin{multline}
-\psi'(\tau)\tau^2\bigg|_{\sigma_{k+2}}^{\sigma_k}+p\psi(\tau)\tau\bigg|_{\sigma_{k+2}}^{\sigma_k}+\int_{\sigma_{k+2}}^{\sigma_k}
\bigg((q-p+\frac{z(\psi(\tau))}{\tau^2})\psi(\tau)\bigg)\,d\tau\\
=\int_{\sigma_{k+2}}^{\sigma_k} \bigg(qs_+-\frac{F'(s_+)}{\tau^2}\psi(\tau)\bigg)\,d\tau.
	\label{eq:psiintegrated}
\end{multline}
Dividing (\ref{eq:psiintegrated}) by $\sigma_k-\sigma_{k+2}$, using (\ref{eq:psisigma}), the existence of $\beta=\lim_{\tau\to 0}\frac{\psi(\tau)}{\tau^2}$ and $\lim_{\tau\to 0} z(\tau)=0$ and then letting $k\to+\infty$, we obtain $qs_+-F'(s_+)\beta=0$.
\end{proof}


\section{Uniqueness and monotonicity}
\label{sec:uniqueness}

\subsection{Uniqueness under positivity assumption}

In our argument, it is more convenient to consider solutions (\ref{RS::ODE})-(\ref{BC}) which satisfy in addition that 
\begin{equation}
u \geq 0 \text{ in }(0,R).
	\label{StripCond}
\end{equation}
See subsection \ref{ssec:NoPosAs} for a discussion on this condition.

The following result gives a statement regarding the range of $u$.
\begin{lemma}\label{Lem:SSC}
Assume \eqref{cond:pq} and \eqref{condF}. If $u$ is a solution of \eqref{RS::ODE}-\eqref{BC} and $u$ satisfies \eqref{StripCond}, then
\begin{equation}
0 < u < s_+ \text{ in }(0,R)
	.\label{StrStripCond}
\end{equation}
\end{lemma}

\begin{proof}  The proof was done in Step $3$ of Proposition~\ref{Prop:OAsymp}.
\end{proof}

A key ingredient in our argument is a comparison principle for the nonlinear ODE \eqref{RS::ODE}. We adopt the following definition for sub/super-solutions of \eqref{RS::ODE}.

\begin{definition}
\label{def1}
A locally Lipschitz, piecewise $C^2$ function $\psi$ defined on a non-empty interval $I$ is said to be a super-solution (or sub-solution) of \eqref{RS::ODE} if it satisfies in $I$
\begin{align*}
&\psi''(r) + \frac{p}{r}\,\psi'(r) - \frac{q}{r^2}\,\psi(r) \leq F(\psi(r)) 
,\\
\Big(\qquad \text{ or \qquad } &\psi''(r) + \frac{p}{r}\,\psi'(r) - \frac{q}{r^2}\,\psi(r) \geq  F(\psi(r)) 
 \qquad \Big)
\end{align*}
wherever it is $C^2$, and if, whenever the first derivative of $\psi$ jumps, says at $r_0 \in I$, there holds
\[
\psi'(r_0^-) >  \psi'(r_0^+) \qquad (\text{or } \psi'(r_0^-) <  \psi'(r_0^+) ).
\]
\end{definition}

We prove: 
\begin{proposition}\label{Prop:CompPrin}
Assume \eqref{cond:pq} and \eqref{condF}. Assume that $\overline{u}$ is a locally Lipschitz, piecewise $C^2$ super-solution of \eqref{RS::ODE} and $\underline{u}$ is a locally Lipschitz, piecewise $C^2$ sub-solution of \eqref{RS::ODE} in $[0,\infty)$. Assume furthermore that
\begin{align*}
&0 \leq \overline{u},\underline{u} \leq s_+,\\
&\overline{u} = \overline{\alpha}\,r^{\ga} + o(r^\ga), \underline{u} = \underline{\alpha}\,r^\ga + o(r^\ga)  \text{ as } r \rightarrow 0,\\
&\overline{u} = s_+ - \overline{\beta}\,r^{-2} + o(r^{-2}), \underline{u} = s_+ - \underline{\beta}\,r^{-2} + o(r^{-2})  \text{ as } r \rightarrow \infty,
\end{align*}
where $\overline{\alpha} > 0$, $\underline{\beta} > 0$ and
\begin{equation}
\underline{\beta} \geq \overline{\beta}.
	\label{Hierarchy}
\end{equation}
Then
\[
\overline{u} \geq \underline{u} \text{ in }(0,\infty).
\]
Moreover, if equality happens somewhere in $(0,\infty)$, then $\overline{u} \equiv \underline{u}$.
\end{proposition}

\begin{proof} {\it Step 1.} We first prove the result under an additional assumption that
\begin{equation}
\overline{\alpha} \geq \underline{\alpha}.
	\label{Eq:alHi}
\end{equation}
We will use the logarithmic sliding method, a variant of the method of moving planes, developed through the works of Alexandrov \cite{Alexandrov}, Serrin \cite{Se}, Gidas, Ni and Nirenberg \cite{G-N-N-1979}, \cite{G-N-N-1981}, and Berestycky and Nirenberg \cite{BN88, BN90}.
Before we begin, we note that, by the argument that led to \eqref{StrStripCond},
\begin{equation}
\overline{u} > 0 \text{ and } \underline{u} < s_+ \text{ in } (0,\infty)
	.\label{StrictBndOUu}
\end{equation}
For any $\theta >0$ we define
\[
\overline{u}_\theta(r) = \overline{u}\Big(\frac{r}{\theta}\Big).
\]
Using $0 \leq \overline{u} \leq s_+$, it is easy to check that, for $\theta < 1$, $\overline{u}_\theta$ is a super-solution to \eqref{RS::ODE}. In fact, by \eqref{condF} and \eqref{StrictBndOUu}, for $\theta < 1$, $\overline{u}_\theta$ is a strict super-solution in the sense that
\begin{equation}
\overline{u}_{\theta}''(r) + \frac{p}{r}\,\overline{u}_{\theta}'(r) - \frac{q}{r^2}\,\overline{u}_{\theta}(r) <  F(\overline{u}_{\theta}(r)) 	\label{StrictSuper}
\end{equation}
wherever $\overline{u}_{\theta}$ is $C^2$. 

\medskip

Our aim is to show that $\overline{u}_\theta \geq \underline{u}$ for any $\theta \in (0,1]$. As consequence, one has $\overline{u} \geq \underline{u}$.

\medskip

\nd {\it Step 1(a). We prove that there exists $\theta_0 > 0$ such that $\overline{u}_\theta > \underline{u}$ in $(0,\infty)$ for any $0 < \theta < \theta_0$.}
By hypotheses, for any $0 < \rho \ll \min(\overline{\alpha},\underline{\beta})$, there exists $\delta_0 = \delta_0(\rho)>0$ such that
\begin{align}
&\overline{u}(r) \geq (\overline{\alpha} - \rho)\,r^\ga
	\text{ and } \underline{u}(r) \leq (\underline{\alpha} + \rho)\,r^\ga
		\text{ for } r < \delta_0
	,\label{bound:uatzero}\\
&\overline{u}(r) \geq s_+ - (\overline{\beta} + \rho)\,r^{-2} 
	\text{ and } \underline{u}(r) \leq s_+ - (\underline{\beta} - \rho)\,r^{-2}
	\textrm{ for }r>\frac{1}{\delta_0}
	.\label{bound:uatinfty}
\end{align}
Replacing $\delta_0$ by some smaller $\tilde\delta_0 < \delta_0$ if necessary, we can further assume that
\begin{equation}\label{ineq:delta0alpharho}
\delta_0^2<\frac 1 4, \quad \max\{\overline{\alpha}-\rho, \underline{\alpha}+\rho\} \delta_0^\ga+(\overline{\beta}+\rho)\delta_0^2\leq s_+.
\end{equation}
From now on, we fix $\rho$ (and so $\delta_0$). For $\delta \in (0,\delta_0]$, define
$$\mathcal{E}(\delta)=\inf_{r\in (\delta,\frac{1}{\delta})} \overline{u}(r).$$

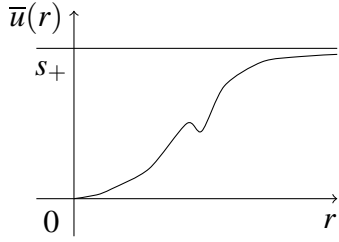
\begin{figure}[h]
\begin{center}
\begin{tikzpicture}
\draw[->] (-0.5,0) -- (3.5,0);
\draw[->] (0,-0.5)--(0,2.5);
\draw (-0.5,2) -- (3.5,2);
\draw plot[smooth] coordinates{(0,0) (0.3,0.05) (0.5,0.13) (1,0.4) (1.5,1) (1.7,0.9) (2,1.5) (2.5,1.82) (3,1.89) (3.5,1.92)};
\draw (3.4,-0.3) node {$r$}
(-0.3,1.7) node {$s_+$}
(-0.3,-0.3) node {$0$}
(-0.5,2.4) node {$\overline{u}(r)$};
\end{tikzpicture}
\end{center}
\caption{A schematic graph of $\overline{u}$.}
\label{fig:fig3}
\end{figure}

Since $\overline{u}$ is locally Lipschitz, \eqref{StrictBndOUu} implies that 
\begin{equation}
\mathcal{E}(\delta) > 0 \text{ for any }\delta \in (0,\delta_0].
\label{OUNontri}
\end{equation}
Using \eqref{bound:uatzero}, \eqref{bound:uatinfty} and \eqref{ineq:delta0alpharho}, one has
\begin{equation}\label{ineq:edelta}
 \mathcal{E}(\delta) \geq \min\Big(\mathcal{E}(\delta_0),(\overline{\alpha}-\rho)\delta^\ga, s_+ - (\overline{\beta} + \rho)\,\delta_0^2\Big) = \min\Big(\mathcal{E}(\delta_0),(\overline{\alpha}-\rho)\delta^\ga\Big),\forall \delta\le \delta_0
	.
\end{equation}

\begin{claim} 
\label{claim1}
Then there exists $\theta_0$ such that, for $0 < \theta < \theta_0$, there holds  $\overline{u}_\theta >  \underline{u}$ in $(0,\infty)$. In fact,
\begin{multline}\label{delta restriction}
\theta_0 := 
	\min\Big(\delta_0^2, 
		\bigg(\frac{\overline{\alpha}-\rho}{\underline{\alpha}+\rho}\bigg)^{1/\ga},
		\bigg(\frac{\overline{\alpha}-\rho}{\underline{\alpha}+\rho}\bigg)^{2/\ga},
		\bigg(\frac{\underline{\beta}-\rho}{\overline{\beta}+\rho}\bigg)^{1/2},\\
		\bigg(\frac{\mathcal{E}(\delta_0)}{\overline{\alpha} - \rho}\bigg)^{2/\ga},
		\bigg(\frac{\mathcal{E}(\delta_0)}{\underline{\alpha} + \rho}\bigg)^{2/\ga},
		\frac{s_+- \sup_{r\in(0,\delta_0^{-1})}\underline{u}(r)}{\overline{\beta}+\rho}\Big).
\end{multline}
(Note that $\theta_0 > 0$ thanks to \eqref{StrictBndOUu}.)
\end{claim}

\begin{proof} Let $\theta = \delta^2$. We check the inequality $\overline{u}_\theta>\underline{u}$ on different intervals:
\begin{itemize}
\item For $r\in (0,\delta^3)$, we have $\frac{r}{\theta} \in (0,\delta)$ and so \eqref{bound:uatzero} and \eqref{delta restriction} give
$$
\overline{u}_\theta(r)
	=\overline{u}\Big(\frac{r}{\theta}\Big)
	\ge (\overline{\alpha}-\rho)\,\frac{r^\ga}{\theta^\ga}
	> (\underline{\alpha}+\rho)r^\ga
	\ge \underline{u}(r).
$$ 

\item For $r\in [\delta^3,\delta^2)$, we have $\frac{r}{\theta}\in [\delta,1)$ and so, by \eqref{bound:uatzero}, \eqref{ineq:edelta}  and \eqref{delta restriction},
\begin{multline*}
\overline{u}_\theta(r)
	= \overline{u}\Big(\frac{r}{\theta}\Big)
	\ge \mathcal{E}(\delta)
	\ge \min\Big(\mathcal{E}(\delta_0),(\overline{\alpha}-\rho)\delta^\ga\Big)
	= (\overline{\alpha} -\rho)\delta^\ga\\
	> (\underline{\alpha}+\rho)\delta^{2\ga}
	\ge (\underline{\alpha} +\rho)r^\ga
	\ge \underline{u}(r)
	.
\end{multline*}

\item For $r\in [\delta^2,\delta)$, we have $\frac{r}{\theta} \in [1,\frac{1}{\delta})$, and so, by \eqref{bound:uatinfty} and \eqref{ineq:delta0alpharho},
$$
\overline{u}_\theta(r) 
	= \overline{u}\Big(\frac{r}{\theta}\Big)
	\ge \min\Big(\mathcal{E}(\delta_0),s_+ - (\overline{\beta} + \rho)\,\delta_0^2\Big)
	> (\underline{\alpha} + \rho)\delta^\ga
	\ge \underline{u}(r).
$$

\item For $r\in [\delta,\frac{1}{\delta_0})$, we have $\frac{r}{\theta} \geq \frac{1}{\delta}$ and so, by \eqref{bound:uatinfty} and \eqref{delta restriction},
$$
\overline{u}_\theta(r)
	=\overline{u}\Big(\frac{r}{\theta}\Big)
	\ge s_+ - (\overline{\beta} + \rho)\frac{\theta^2}{r^2}
	\geq s_+ - (\overline{\beta} + \rho)\delta^2
	> \underline{u}(r).
$$

\item Finally, for $r\in (\frac{1}{\rho_0},\infty)$ we have $\frac{r}{\theta} > \frac{1}{\theta}$ and so, by \eqref{bound:uatinfty}  and \eqref{delta restriction},
$$
\overline{u}_\theta(r)
	=\overline{u}\Big(\frac{r}{\theta}\Big)
	\ge s_+- (\overline{\beta}+\rho)\frac{\theta^2}{r^2}
	> s_+- \frac{\underline{\beta}-\rho}{r^2}
	\ge \underline{u}(r).
$$

\end{itemize}

\nd We have thus shown that $\overline{u}_\theta > \underline{u}$ for any $\theta \in (0,\theta_0)$ which ends the proof of Step 1. 
\end{proof}

\medskip

\nd {\it Step 1(b).}
Define 
\[
\bar\theta = \sup\{\theta < 1: \overline{u}_\sigma \geq \underline{u} \text{ in } (0,\infty),\forall 0 < \sigma \leq \theta\}
	\;.
\]
Evidently, $\bar\theta$ is well-defined, $\theta_0 \leq \bar \theta \leq 1$ and $\overline{u}_{\bar\theta}\geq \underline{u}$ in $(0, \infty)$. To complete the proof,  we need to show that $\bar \theta =1$. 

\begin{claim}
\label{claim2}
If $\bar\theta <1$, then there exists $r_0\in (0, \infty)$ such that $\overline{u}_{\bar\theta}(r_0) = \underline{u}(r_0)$.
\end{claim}

\begin{proof} 
Arguing indirectly, assume that $\overline{u}_{\bar\theta} > \underline{u}$ on $(0,\infty)$. To get a contradiction, we show that there exists $\mu_0 > 0$ such that $\overline{u}_{\bar\theta + \mu} \geq \underline{u}$ for any $0 < \mu < \mu_0$. Select $\eps > 0$ and $0 < \mu_1 < 1 - \bar\theta$ such that
\be
\label{num}
\frac{\overline{\alpha} - \eps}{(\bar\theta + \mu_1)^\ga} > \underline{\alpha} + \eps \text{ and } (\overline{\beta} + \eps)(\bar\theta + \mu_1)^2 < \underline{\beta} - \eps.
\ee
Such $\eps$ exists thanks to \eqref{Hierarchy} and that $\bar\theta < 1$.
By \eqref{bound:uatzero}, we have for $0 < r < \bar\theta\,\delta_0(\eps)$ and $0 < \mu < \mu_1$ that
\[
\overline{u}_{\bar\theta + \mu}(r)
	= \overline{u} \Big(\frac{r}{\bar\theta + \mu}\Big)
	\geq (\overline{\alpha} - \eps)\frac{r^\ga}{(\bar\theta + \mu)^\ga}
	\stackrel{\eqref{num}}{>} (\underline{\alpha} + \eps)r^\ga
	\geq \underline{u}(r).
\]
Likewise, by \eqref{bound:uatinfty}, we have for $r > \frac{1}{\delta_0(\eps)}$ and $0 < \mu < \mu_1$ that
\[
\overline{u}_{\bar\theta + \mu}(r)
	= \overline{u} \Big(\frac{r}{\bar\theta + \mu}\Big)
	\geq s_+ - (\overline{\beta} + \eps)\,\frac{(\bar\theta + \mu)^2}{r^2}
	\stackrel{\eqref{num}}{>} s_+ - (\underline{\beta} - \eps)\frac{1}{r^2}
	\geq \underline{u}(r).
\]
On the other hand, since $\overline{u}_{\bar\theta} > \underline{u}$ in $[\bar\theta \delta_0(\eps), \frac{1}{\delta_0(\eps)}]$, which is compact, we can select $\mu_2 > 0$ sufficiently small such that for any $0 < \mu < \mu_2$, there holds that $\overline{u}_{\bar\theta + \mu} > \underline{u}$ in $[\bar\theta \delta_0(\eps), \frac{1}{\delta_0(\eps)}]$. Altogether, we just showed that if $0 < \mu < \mu_0 = \min(\mu_1,\mu_2)$, then 
\[
\overline{u}_{\bar\theta + \mu} > \underline{u} \text{ in } (0,\infty).
\]
This contradicts the maximality of $\bar\theta$. Therefore, there exists $r_0\in (0, \infty)$ such that $\overline{u}_{\bar\theta}(r_0) = \underline{u}(r_0)$.
\end{proof}

\begin{claim}
\label{claim3}
If $\bar\theta \leq1$ and there exists $r_0\in (0, \infty)$ such that $\overline{u}_{\bar\theta}(r_0) = \underline{u}(r_0)$, then $\bar\theta=1$.
\end{claim}

\begin{proof}
Recalling the definition of super/sub-solutions, the equality $\overline{u}_{\bar\theta}(r_0) = \underline{u}(r_0)$ forces the first derivatives of $\overline{u}_{\bar\theta}$ and 
$\underline{u}$ to be continuous across $r_0$. Consider the function $w=\overline{u}_{\bar\theta} - \underline{u}$. Then $w$ has a local minimum at $r_0$, is $C^1$-continuous at $r_0$ and possesses left and right second derivatives at $r_0$. In addition, as $\overline{u}_{\bar\theta}$ is a super-solution while $\underline{u}$ is a sub-solution, we deduce that $w''(r_0^\pm) \leq 0$. This forces $w''(r_0^\pm) = 0$. Hence $w$ is $C^2$ across $r_0$ and so in a neighborhood, say $(r_-,r_+)$, of $r_0$. Observe that $w$ satisfies
$$
w'' + {p \over r} w' - {q \over r^2} w \leq c(x) w , \qquad  w \geq 0 \text{ in } (r_-,r_+) \text{ and } w(r_0) = 0,
$$
where $c(x) = { F(\overline{u}_{\bar\theta}) - F(\underline{u}) \over {\overline{u}_{\bar\theta} -\underline{u}}}(x)$. The strong maximum principle
then implies that $w \equiv 0$ in $(r_-,r_+)$. In other words, $\overline{u}_{\bar\theta} \equiv \underline{u}$ in $(r_-,r_+)$. It is readily seen that this statement implies that $\overline{u}_{\bar\theta} \equiv \underline{u}$ in $(0,\infty)$. In particular, $\overline{u}_{\bar\theta}$ is a solution of \eqref{RS::ODE} in $(0,\infty)$. Recalling \eqref{StrictSuper}, it follows that $\bar\theta = 1$. This ends the proof of Claim \ref{claim3}.
\end{proof}

\nd By Claims \ref{claim2} and \ref{claim3}, we deduce that $\overline{u} \geq \underline{u}$ in $(0,\infty)$. The rigidity statement follows from the proof of Claim \ref{claim3}. We have thus proved the assertion when \eqref{Eq:alHi} holds.

\medskip
\noindent{\it Step 2.} To complete the proof, we prove \eqref{Eq:alHi}. Assume by contradiction that $\overline{\alpha} < \underline{\alpha}$. Define $\overline{u}_\theta(r) = u(r/\theta)$ as above. We have seen that, for $0 < \theta \leq 1$, $\overline{u}_\theta$ is a super-solution.

Select $\theta$ such that $\overline{\alpha}\theta^{\gamma_+} = \underline{\alpha}$. Applying the result obtained in Step 1 for $\tilde u:= \overline{u}_\theta$ and $\underline{u}$, we obtain $\tilde u \geq \underline{u}$.

Let $v = \tilde u - \underline{u} \geq 0$. Then $v$ satisfies
\[
v'' + \frac{p}{r} v' - \frac{q}{r^2} v + c(r)v \leq 0 \text{ and } \lim_{r \rightarrow 0} \frac{v}{r^{\gamma_+}} = 0,
\]
where $c$ is some function which is continuous in $[0,\infty)$. Let $w =  \frac{v}{r^{\gamma_+}}$, then $w$ satisfies
\[
w'' + \frac{p + 2\gamma_+}{r} w'  + c(r)w \leq 0 \text{ and } w(0) = 0.
\]
Since $w \geq 0$ and $p + 2\gamma_+ > 1$, Lemma \ref{Lem:SingHopf} implies that $w \equiv 0$, i.e. $\tilde u \equiv \underline{u}$. This forces $\theta = 1$ and so $\overline u \equiv \tilde u \equiv \underline{u}$, which contradicts the assumption that $\overline{\alpha} < \underline{\alpha}$. We have thus proved \eqref{Eq:alHi} and completed the proof of the proposition.
\end{proof}

\begin{remark}\label{RemW}
The conclusion of Proposition \ref{Prop:CompPrin} remains valid if one replaces the condition
\[
\underline{u} = s_+ - \underline{\beta}\,r^{-2} + o(r^{-2})  \text{ as } r \rightarrow \infty
\]
by the condition
\[
\limsup_{r \rightarrow \infty} \underline{u} < s_+.
\]
\end{remark}

As a consequence of the argument in Step 2 of the proof, we have the following Hopf-type lemma:

\begin{corollary}\label{ERS::Hopf}
Assume that $\overline{u}$ is a super-solution of \eqref{RS::ODE} and $\underline{u}$ is a sub-solution of \eqref{RS::ODE} in $[0,R)$ for some $0 < R < \infty$ such that both can be factored as a product of $r^\ga$ and a continuous function at $r = 0$. If $\overline{u} \geq \underline{u}$ in $(0,R)$ then
\[
\text{ either } \lim_{r\rightarrow 0} \frac{\overline{u}}{r^\ga} > \lim_{r\rightarrow 0} \frac{\underline{u}}{r^\ga} \text{ or } \overline{u} \equiv \underline{u}.
\]
\end{corollary}

The following results are variants of the previous comparison principle on different intervals.

\begin{proposition}\label{Prop:CompPrinFinite}
Assume \eqref{cond:pq} and \eqref{condF}. Assume that $\overline{u}$ is a locally Lipschitz, piecewise $C^2$ super-solution of \eqref{RS::ODE} and $\underline{u}$ is a locally Lipschitz, piecewise $C^2$ sub-solution of \eqref{RS::ODE} on some interval $\mathcal{I}\subset (0,\infty)$.

{\it (i)}  Assume that $\mathcal{I}=(0,R)$ with $R<\infty$. Furthermore assume that
\begin{align*}
&0 \leq \overline{u},\underline{u} \leq s_+,\\
&\overline{u} = \overline{\alpha}\,r^\ga + o(r^\ga), \underline{u} = \underline{\alpha}\,r^\ga + o(r^\ga)  \text{ as } r \rightarrow 0,\\
&\overline{u}(R) = s_+,  \underline{u}(R) \leq s_+,
\end{align*}
where $\overline{\alpha} > 0$. Then
\[
\overline{u} \geq \underline{u} \text{ in }(0,R).
\]

\medskip
{\it (ii)} Assume that $\mathcal{I}=(r_1,\infty)$ with $0 \leq r_1 < \infty$. Furthermore assume that
\begin{align*}
&0 \leq \overline{u},\underline{u} \leq s_+,\\
&\overline{u}(r_1) \geq 0, \underline{u}(r_1) = 0,\\
&\overline{u} = s_+ - \overline{\beta}\,r^{-2} + o(r^{-2}), \underline{u} = s_+ - \underline{\beta}\,r^{-2} + o(r^{-2})  \text{ as } r \rightarrow \infty,
\end{align*}
where $\underline{\beta}$ and $\overline{\beta}$ satisfy \eqref{Hierarchy}.
Then
\[
\overline{u} \geq \underline{u} \text{ in }(r_1,\infty).
\]

Moreover, in either case we have that if equality happens somewhere in $\mathcal{I}$, then $\overline{u} \equiv \underline{u}$.
\end{proposition}

\begin{proof}
{\it (i)} The proof goes exactly the same, but simpler, as in that of Proposition \ref{Prop:CompPrin}. The key difference is that $\overline{u}_\theta$ ($0 < \theta < 1$) is defined by
\[
\overline{u}_\theta(r) = \left\{\begin{array}{ll}
	\overline{u}\big(\frac{r}{\theta}\big) \text{ for } 0 \leq r < \theta\,R,\\
	s_+ \text{ for } \theta\,R \leq r \leq R.
\end{array}\right.
\]
We omit the details.

\medskip
{\it (ii)} Again the proof is a variant of that of Proposition \ref{Prop:CompPrin}. First extend $\underline{u}$ by setting
\[
\underline{u}(r) = 0 \text{ for } 0 < r < r_1.
\]
Note that the extended function $\underline{u}$ is a sub-solution of \ref{RS::ODE} on the whole interval $(0,\infty)$. Next, define
\[
\overline{u}_\theta(r) =
	\overline{u}\big(\frac{r}{\theta}\big) \text{ for } \theta\,r_1 \leq r < \infty.
\]
Then $\overline{u}_\theta$ is a super-solution of \eqref{RS::ODE} in $(\theta\,r_1, \infty)$ for all $\theta \in (0,1)$. The proof of Proposition \ref{Prop:CompPrin} can now be applied to reach the conclusion. We omit the details.
\end{proof}

Combining Propositions \ref{Prop:OAsymp}, \ref{ERS::Asymp}, \ref{Prop:CompPrin} and \ref{Prop:CompPrinFinite}  we obtain the following uniqueness statements.

\begin{proposition}\label{ERS::Uniqueness}
Assume that $p$ and $q$ satisfies \eqref{cond:pq} and $F$ satisfies \eqref{condF}. For any $0 < R \leq \infty$, there is at most one non-negative solution $u$ to the BVP \eqref{RS::ODE}\&\eqref{BC}.
\end{proposition}

To conclude the section, we turn to monotonicity properties for solutions of \eqref{RS::ODE}\&\eqref{BC}.
\begin{lemma}\label{ERS:monotonicity}
For any $0 < R \leq \infty$, if $u$ is a solution of (\ref{RS::ODE})\&(\ref{BC}), and $r_1 \in [0,R)$ is the last zero of $u$ (i.e. $u(r_1) = 0$ and $u(r) > 0$ for $r \in (r_1,R)$), then $u$ is strictly increasing in $(r_1,R)$.
\end{lemma}

\begin{proof} Let us consider first the case when $r_1=0$ and $R = \infty$. By Proposition \ref{Prop:OAsymp}, $u$ can be expressed as a product of $r^\ga$ and a continuous function at $r = 0$. Recalling \eqref{StrStripCond}, we can apply Corollary \ref{ERS::Hopf} to obtain
\[
\lim_{r\rightarrow 0}\frac{u}{r^\ga} > 0.
\]

Now, for any $\theta >0$ we define
\[
u_\theta(r) = u\Big(\frac{r}{\theta}\Big)
	\;.
\]
Using \eqref{StrStripCond}, it is easy to check that $u_\theta$ is a super-solution of \eqref{RS::ODE} for $0 < \theta < 1$. Keeping in mind Proposition \ref{ERS::Asymp}, we can apply the comparison principle in Proposition \ref{Prop:CompPrin} to $\overline{u} = u_\theta$ and $\underline{u} = u$ to conclude that
\[
u_{\theta}(r) > u(r) \text{ for any } 0 < r < \infty \text{ and } 0  < \theta < 1.
\]
In particular, for $0 < r < s < \infty$,
\[
u(r) < u_{\frac{r}{s}}(r) = u(s).
\]
This completes the proof for the case $R = \infty$.

The proof in the case $r_1 = 0$ and $R < \infty$ is similar:  One applies the comparison principle in Proposition \ref{Prop:CompPrinFinite} to $\underline{u} = u$ and $\overline{u} = u_{\theta}$ where this time $u_{\theta}$ is defined by
\[
u_\theta(r) = \left\{\begin{array}{ll}
	u\big(\frac{r}{\theta}\big) \text{ for } 0 \leq r < \theta\,R,\\
	s_+ \text{ for } \theta\,R \leq r \leq R.
\end{array}\right.
\]
We omit the details.

\medskip Assume now that $r_1 > 0$. We present the proof for the case $R = \infty$. The case $R < \infty$ can be done similarly.

For any $\theta \in (0,1)$ we define
\[
u_\theta(r) = u\Big(\frac{r}{\theta}\Big) \text{ for } r \geq \frac{r_1}{\theta}.
	\;.
\]
Then $u_\theta$ is a super-solution of \eqref{RS::ODE} in $(r_1/\theta, \infty)$. On the other hand, if we set
\[
\underline{u}(r) = \left\{\begin{array}{l}
0 \text{ for } r \in (0,r_1),\\
u(r) \text{ for } r \in [r_1,\infty),
\end{array}\right.
\]
then $\underline{u}$ is a sub-solution of \eqref{RS::ODE} in $(0,\infty)$. We can then apply Proposition \ref{ERS::Asymp} and the comparison principle in Proposition \ref{Prop:CompPrinFinite} to $\overline{u} = u_\theta$ and $\underline{u}$ to conclude the proof.
\end{proof}

\bigskip We can now gather previously developed arguments to present:

\begin{proof}[Proof of Theorem~\ref{thm:main}.]
The existence of the solution for the case $R<\infty$ is a consequence of Lemma~\ref{lemma:uR}, where the solution is obtained as a global energy minimizer of the modified energy $\tilde E$ defined in \eqref{energyODE0}. In Corollary~\ref{cor:energyuniqueness} next section it will be noted  
that if the nonlinearity $F$ satisfies condition  \eqref{Fcond:even}  then the solutions thus obtained are global energy minimizers of the standard energy \eqref{energyODE}.

In the case of infinite domain, $R=\infty$, the existence of the solutions is obtained in Proposition~\ref{prop:existenceinftydomain} as limit of solutions obtained for finite $R$, as $R\to\infty$. The most delicate part is to ensure that the solution thus obtained satisfies the boundary conditions at $0$ and $\infty$. In order to study the behaviour at $0$ we use Proposition~\ref{Prop:OAsymp}, while in order to study the asymptotics at $\infty$ we use the monotonicity results  of Lemma~\ref{ERS:monotonicity} together with an energy argument which also shows that the the solution thus obtained is locally energy minimizing.

In order to prove uniqueness we first show in Lemma~\ref{Lem:SSC} which provides that a non-negative solution is actually positive and stays away from $s_+$, and use these in the study of sub-solutions and super-solutions in  Lemma~\ref{Prop:CompPrin}. Combining this last lemma with the detailed behaviour at $0$ obtained in Proposition~\ref{Prop:OAsymp} and the one at $\infty$ obtained in Proposition~\ref{ERS::Asymp} we obtaine the uniqueness of positive solutions stated in Proposition~\ref{ERS::Uniqueness}.
\end{proof}


\subsection{Uniqueness without positivity assumption}

\label{ssec:NoPosAs}
  
   In this section we consider  two different  types of additional assumptions under which we can obtain the uniqueness of solutions for  \eqref{RS::ODE} and \eqref{BC}, without the positivity requirement on $u$.

 The first condition is imposed on the solution while the second one is a condition on the nonlinearity. In either case we  show that in fact a nodal solution must necessarily be positive and then the uniqueness result in class of positive solutions will provide us the more general uniqueness result.

\smallskip
We start by noting that positivity is implied by the requirement of local energy minimization, as stated in  Corollary~\ref{cor:energyuniqueness}. We now show:

\medskip
\begin{proof}[Proof of Corollary~\ref{cor:energyuniqueness}.]
We claim that assumption \eqref{Fcond:even} implies that the solution $u$ obtained in Theorem~\ref{thm:main} is locally energy minimizing. Indeed, as in Section~\ref{ssec:ExistenceFinite}, since $F$ satisfies \eqref{FCond:even}, then Lemma~\ref{lemma:uR} provides the claim in the case of bounded domains $(0, R)$, $R<\infty$. In the case of unbounded domain, the solution $u_\infty$ obtained in the proof of Proposition~\ref{prop:existenceinftydomain} is locally energy minimizing.

\medskip
We consider now the converse: we take $u \in H^1_{\loc}(0,R)$ that is a locally energy minimizing solution of \eqref{RS::ODE} with respect to the  energy \eqref{energyODE} and satisfies $u(R)=s_+$.

Arguing similarly as in the proof of Lemma~\ref{lemma:uR} we have $u$ and $|u|$ are both minimizers for $E$ on $(0,R')$ for all sufficiently large $R' < R$ such that $u(R') > 0$. Thus $|u|$ is a non-negative solution of \eqref{RS::ODE}. As shown in Step 3 of the proof of Proposition \ref{Prop:OAsymp}, this implies that $|u| > 0$ in $(0,\infty)$, and, as $u$ has constant sign, we have $u > 0$ in $(0,\infty)$.

We also claim that $u < s_+$ in $(0,R)$. Indeed, if $u(R') \geq s_+$ for some $0 < R' < R$, let
\[
\tilde u(r) = \min (u(r), u(R')), \,0<r<R'
\]
As in the proof of Lemma \ref{lemma:uR}, the fact that $F > 0$ in $(s_+,\infty)$ implies that $E[\tilde u;(0,R')] \leq E[u;(0,R')]$ where equality holds if andy only if $u \equiv u(R')$ in $(0,R')$. Since $u$ is minimizing in $(0,R')$ this implies that $u \equiv u(R')$ in $(0,R')$ which is impossible in view of equation \eqref{RS::ODE}. The claim is proved.

We have proved that $0 < u < s_+$ in $(0,R)$. Since $u$ is bounded, then Proposition \ref{Prop:OAsymp} implies $u(0)=0$. The uniqueness part in Theorem~\ref{thm:main} now shows that $u$ in fact coincides with the solution of \eqref{RS::ODE}\&\eqref{BC} obtained therein. 
\end{proof}

Moving on to imposing conditions on the nonlinearity, we note first that a simple condition  on the behaviour of the nonlinearity on  $(-\infty,0]$ allows to deduce the positivity of any solution of \eqref{RS::ODE} and \eqref{BC}. Indeed, 
if $F$ satisfies $F(t) < 0$ for $t < 0$, then \eqref{StripCond} can be proved using the maximum principle. However this is not satisfied for the physical potential  $F$ of the form \eqref{def:physical nonlinearity}.
To obtain conditions for showing the positivity of solutions for physical type of nonlinearities $F$, we need to impose more constraints on $F$, namely the ones in \eqref{condFLeft}. We now show:

\begin{proof}[Proof of Theorem~\ref{thm:uniqueness+}.]
 We present an argument that is reminiscent of the one in Proposition $3$ in \cite{Hervex2} . For simplicity, we will only present a proof when $R = \infty$. The other case requires only minor modifications. By Step 3 in the proof of Proposition \ref{Prop:OAsymp}, we show that $u < s_+$ in $(0,\infty)$. Using the first line of \eqref{condFLeft}, the same argument shows that $u > s_-$ in $(0,\infty)$. We claim that $u>0$ on $(0,\infty)$ and therefore, by Theorem \ref{thm:main}, $u$ is unique. Arguing by contradiction let us  assume that $u$ is negative somewhere. Since $u(r) \rightarrow s_+$ as $r \rightarrow \infty$, there is some $r_1 \in (0,\infty)$ such that 
\[
u(r_1) = 0 \, \, \text{ and } \, \, u(r) > 0 \text{ for } r > r_1.
\]
In particular, $u'' + \frac{p}{r} u' - \frac{q}{r^2} u  = F(u(r)) < 0$ in $(r_1,\infty)$. By the Hopf lemma, we have
\begin{equation}
	\label{Eq:dur1>0}
u'(r_1) > 0.
\end{equation}
Hence there exists $r_0 \in [0,r_1)$ such that $u(r_0) = 0$ and $u(r) < 0$ for $r \in (r_0,r_1)$.

\smallskip
We now define $r_2 = 2r_1 - r_0$ and $\psi(r):=-u(2r_1-r)$ for $r\in (r_1, r_2)$. Then $\psi$ is positive in $(r_1,r_2)$, $\psi(r_1) = \psi(r_2) = 0$ and $\psi$ satisfies the ODE:
\[
\psi''-\frac{p}{2r_1-r}\psi'-\frac{q}{(2r_1-r)^2}\psi=-F(-\psi), \quad r\in (r_1, r_2).
\]
In addition,
\[
\psi(r_1)=u(r_1)=0,\,\psi'(r_1)=u'(r_1)>0,\,\psi''(r_1)=-u''(r_1)=\frac{p}{r_1} u'(r_1)>0.
\]
Thus, for some $\varepsilon \in (0, r_1-r_0)$, we have $\psi >u$ on $(r_1,r_1 + \varepsilon)$. Let $r_3\in (r_1, r_2)$ be the maximal point where $\psi >u$ on $(r_1,r_3)$, so that
$\psi(r_3)=u(r_3)$ (this is possible because $\psi(r_2)=0<u(r_2)$).
On $(r_1,r_3)$ we have
\begin{align}
(u'\psi-u\psi')'=&u\psi \left[\underbrace{\frac{q}{r^2}-\frac{q}{(2r_1-r)^2}}_{<0\,\textrm{ for }r>r_1}+\underbrace{\frac{F(u)}{u}+\frac{F(-\psi)}{\psi}}_{\le 0\textrm{ by }\eqref{condFLeft}}\right]-\frac pr u'\psi-\frac{p}{(2r_1-r)}\psi' u\nonumber\\
&\le \frac{p}{2r_1-r}\left(u'\psi-u\psi'\right)
\label{inegu}
\end{align} where for the last inequality we used that $u'>0$ (see Lemma~\ref{ERS:monotonicity}), $\psi, u\ge 0$ on $(r_1, r_3)$ and $p>0$.
If we denote $\zeta(r):=u'(r)\psi(r)-u(r)\psi'(r)$ and $f(r):=\frac{p}{2r_1-r}$ then \eqref{inegu} implies $\zeta'(r)\le f(r) \zeta(r)$ on $(r_1, r_3)$. Noting that $f$ is integrable on $(r_1, r_3)$ and $\zeta(r_1)=0$ we have by Gronwall's inequality that $\zeta\le 0$ on $(r_1,r_3)$. We obtain thus that $\frac{u}{\psi}$ is non-increasing on $(r_1,r_3)$. This leads to a contradiction since $\frac{u}{\psi} < 1$ in $(r_1, r_3)$ while $u(r_3)=\psi(r_3)>0$. 
\end{proof}

\begin{remark}\label{rmk:p01}
Let us point out the modifications needed in the previous argument if $p=0$. We now assume that $F$ is of class $C^2$ and we note that this together with \eqref{condFLeft}  implies $\frac{F(t)+F(-t)-2F(0)}{t^2}\le 0$ for $t$ small enough, hence $F''(0)\le 0$. Then following the previous  proof we first note that  $\psi''(r_1)=-u''(r_1)=0$. In order to compare the behaviours of $\psi$ and $u$ at $0$, we need to compute higher order derivatives as $0$. We have $\psi'''(r_1)=u'''(r_1)=0$ and $\psi^{(4)}(r_1)=-u^{(4)}(r_1)=-\left[ -\frac{4q}{r_1^3} u'(r_1)+F''(0)|u'|^2\right]>0$. The proof continues similarily as before.

\end{remark}

\begin{remark}\label{rmk:p02}
We point out an alternative approach for dealing with the case $p=0$ under a different assumption on the nonlinearity. Namely in addition to \eqref{condF} we require that there exists $\alpha>1$ such that \eqref{condF+} holds. Then in the proof of Theorem~\ref{cor:energyuniqueness} we take a different definition of $\psi$ namely  $\psi(r):=-u\left((\alpha+1)r_1-\alpha r\right)$. We denote $r_2:=\frac{1}{\alpha}\left[(\alpha+1)r_1-r_0\right]$ and observe that $\psi(r_2)=0$. We obtain that $\psi$ satisfies the equation:
\[
\psi''-\frac{\alpha p\psi'}{((\alpha+1)r_1-\alpha r)}-\frac{q\alpha^2\psi}{((\alpha+1)r_1-\alpha r)^2}=-\alpha^2 F(-\psi)
\] and $\psi(r_1)=u(r_1)$, $\psi'(r_1)=\alpha u'(r_1)>0$ hence $\psi'(r_1)>u'(r_1)$ and thus $\psi>u$ on some maximal interval $(r_1,r_3)$. Moreover we have
$$(u'\psi-u\psi')'=u\psi \left[\underbrace{\frac{q}{r^2}-\frac{\alpha^2 q}{((\alpha+1)r_1-\alpha r)^2}}_{\le 0 \text{ for }r\in [r_1,r_2]}+\underbrace{\frac{F(u)}{u}+\frac{\alpha^2 F(-\psi)}{\psi}}_{\le 0\text{ by }\eqref{condF+}}\right]\le 0\quad \textrm{on} \, (r_1, r_3)$$ and this shows that $u/\psi$ is non-increasing of $(r_1,r_3)$. We reach thus a contradiction because $u/\psi<1$ on $(r_1, r_3)$ and $u(r_3)/\psi(r_3)=1$.
\end{remark}

\begin{figure}[h]
\begin{center}
\includegraphics[scale=0.3]{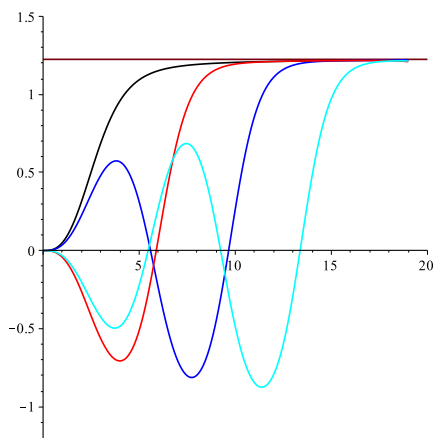}
\end{center}
\caption{A plot of several different solutions for $F(u)=-u+\frac{2}{3}u^3$ and $p=-1,q=3$.}
\label{fig:multiplesol}
\end{figure}

\begin{remark}\label{Rem:pNeg}
If $p<0$ numerical explorations show that there can be several sign changing solutions. See Figure~\ref{fig:multiplesol}.
\end{remark}



\section{Refined qualitative analysis}
\label{sec:refined}
In this section we prove several refined qualitative properties of the positive solution in the physically motivated case. Throughout this section we assume $r \in [0,\infty)$, $p = 2$, $q = 6$, $F$ will take the form \eqref{def:physical nonlinearity} and we denote
\be
u(r) \hbox{ is the unique solution of \eqref{RS::ODE}-\eqref{BC}  and } w(r) = \frac{r u'(r)}{u(r)}.
 \label{not}
\ee

Define
\begin{align*}
f(u) &:= \frac{F(u)}{u} = -a^2 - \frac{b^2}{3}u + \frac{2c^2}{3}u^2,\\
\hat f(u) &:= f'(u)u + f(u) = -a^2 - \frac{2b^2}{3}u + 2c^2\,u^2.
\end{align*}
Then
\begin{align}
u'' + \frac{2}{r} u' - \frac{6}{r^2}u &= uf(u), \label{u-Equation}\\
u''' + \frac{2}{r}u'' - \frac{8}{r^2}\,u' + \frac{12}{r^3}u &= \hat f(u)\,u'.\label{du-Equation}
\end{align}

Note that, by Proposition \ref{Prop:OAsymp}, $v(r)=\frac{u(r)}{r^2}$ is decreasing on the interval $r \in (0, \infty)$ and as a consequence the function $w(r)$ satisfies  $0<w(r) < 2$ for all $r \in (0,\infty)$.

\begin{lemma}
For the function $w(r)$ defined in \eqref{not} the following inequalities hold
\[
2w(w - 2) < rw' < w(w-2) < 0 \text{ in } (0,\infty).
\]
In particular, $w(r)$ and $\frac{2 - w}{w r^4}$ are strictly decreasing and $\frac{2-w}{wr^2}$ is strictly increasing on $(0,\infty)$.
\end{lemma}
\begin{proof} We first show that $w$ is decreasing. 
Straightforward calculations using \eqref{u-Equation} and \eqref{du-Equation} give
\begin{align}\label{eq:A3-1}
w'
	&
	=-\frac{1}{r}(w - 2)(w + 3) + r\,f(u)
	,\\ \label{eq:A3-2}
w''	&= \frac{2}{r^2}(w - 2)(w + 1)(w + 3)  + (\hat f (u) - 3 f(u))\,w
	,\\ \label{eq:A3-3}
w'''
	&=  \frac{w'}{w}\,w'' +  \frac{4}{r^2} \frac{w^3 + w^2 + 3}{w}\,w' - \frac{4}{r^3}(w - 2)(w + 1)(w + 3) + \frac{b^2}{3}u'\,w
	.
\end{align}
Since $0 < w < 2$ we see that on $(0,\infty)$ the function $w'(r)$  satisfies
\be \label{eqw}
(w')'' - \frak{p}(r)\,(w')' - \frak{q}(r)\,w' \geq - \frac{4}{r^3}(w - 2)(w + 1)(w + 3) > 0
\ee
where $\frak{p}(r) = \frac{w'}{w}$ and $\frak{q}(r) = \frac{4}{r^2} \frac{w^3 + w^2 + 3}{w}  > 0$.

By Proposition~\ref{Prop:OAsymp}, $\frac{u}{r^2}$ is decreasing and differentiable up to $r=0$ and its derivative at $0$ is $0$. Thus
$$
\lim_{r \to 0} \frac{u(r)}{r^2} = \alpha > 0 \text{ and } 0 = \lim_{r\to 0} \left[\frac{u(r)}{r^2}\right]'=\lim_{r\to 0}\frac{u(r)}{r^2}\frac{(w(r)-2)}{r}.
$$
It follows that 
\begin{equation}
\lim_{r \rightarrow 0} w(r) = 2  \text{ and } \lim_{r\to 0} w'(r)=0.
	\label{15V13-X1}
\end{equation}
Also, using definition of $w$ and \eqref{asympt:u} we derive that 
\[
\lim_{r\rightarrow \infty} w'(r) = 0.
\]
We can then apply the maximum principle to \eqref{eqw} to conclude that $w' < 0$ on $(0, \infty)$.

We now want to show the upper bound for $rw'(r)$
\[
\chi(r)  := rw' - w(w-2) < 0.
\]
The idea of the proof is essentially the same as before: we find the differential inequality for $\chi(r)$ and employ the maximum principle. A calculation gives
\begin{align*}
\chi'
	&= rw'' + 3w' - 2w\,w'
	,\\
\chi''
	&= rw''' + 4w'' - 2w\,w'' - 2|w'|^2\\
	&\geq (r\frak{p} + 2(2 - w))w'' + r\frak{q}\,w' -   2|w'|^2 -  \frac{4}{r^2}(w - 2)(w + 1)(w + 3)\\
	&= \frac{r\frak{p} + 2(2 - w)}{r}\chi' + \Big[- \frac{3w'}{r w}  +  \frac{3(5w^2 - 2w + 4)}{r^2w}\Big]\chi 
				 +  \frac{11}{r^2}w(w - 2)^2,
\end{align*}
where in the first inequality we have used \eqref{eqw}. Recalling that $w' < 0$ and $0 < w < 2$ on $(0,\infty)$, we see that $\chi$ satisfies
\[
\chi'' - \tilde{\frak{p}}\,\chi' - \tilde{\frak{q}}\,\chi \geq \frac{11}{r^2}w(w - 2)^2 > 0
\]
where $\tilde{ \frak{p}} = \frac{r\frak{p} + 2(2 - w)}{r}$ and $\tilde {\frak{q}} = - \frac{3w'}{r w}  +  \frac{3(5w^2 - 2w + 4)}{r^2w} \geq 0$. In addition, by \eqref{asympt:u}, \eqref{15V13-X1} and the expression for $w'$ we have $\lim_{r \rightarrow 0} \chi (r) = 0 = \lim_{r \to \infty} \chi(r) =0$. Applying the maximum principle we obtain $\chi (r) < 0$ on $(0,\infty)$.

Finally, we show that
\[
\hat \chi := rw' - 2w(w - 2) > 0.
\]
We compute 
\begin{align*}
\hat \chi'
	&= rw'' + 5w' - 4w\,w'
	,\\
\hat \chi''
	&= (r\frak{p} + 2(3 - 2w))w'' + r\frak{q}\,w' -   4|w'|^2 -  \frac{4}{r^2}(w - 2)(w + 1)(w + 3) + \frac{b^2}{3}r\,u'\,w.
\end{align*}
Recalling the definitions of $f(u)$ and $\hat f(u)$ and equation for $w''$ we have
\[
\frac{b^2}{3}\,r\,u'\,w = \frac{b^2}{3}\,u\,w^2 \leq [\hat f(u) - 3f(u)]\,w^2 = w\,w'' - \frac{2}{r^2}w(w - 2)(w + 1)(w + 3).
\]
Using the above inequality and combining the terms we obtain 
\begin{align*}
\hat \chi''
	&\leq (rp + 3(2 - w))w'' + rq\,w' -   4|w'|^2 -  \frac{2}{r^2}(w - 2)(w + 1)(w + 3)(w + 2)\\
	&= \frac{r\frak{p} + 3(2 - w)}{r}\hat \chi' - \Big[\frac{5w'}{rw} + \frac{8w^3 - 33 w^2 + 10w - 12}{r^2w}\Big]\hat \chi\\
		&\qquad   - \frac{6}{r^2}(w - 2)^2(3w^2 - 3w + 1)
		.
\end{align*}
Recalling that $0 < w < 2$ and an upper bound for $w'$ we see that
\[
\hat \chi'' - \hat{\frak{p}}\,\hat\chi' - \hat{\frak{q}}\,\hat\chi \leq - \frac{6}{r^2}(w - 2)^2(3w^2 - 3w + 1)  < 0
\]
where $\hat{\frak{p}}= \frac{rp + 2(3 - 2w)}{r}$ and $\hat{\frak{q}} = - \frac{5w'}{rw} - \frac{8w^3 - 33 w^2 + 10w - 12}{r^2w} > 0$. As in the  previous case we also have $\lim_{r \rightarrow 0} \hat \chi(r) = 0 = \lim_{r \to \infty} \hat \chi(r)$ and so the maximum principle gives $\hat \chi >0$ on $(0, \infty)$.
\end{proof}

\begin{remark}\label{remark:fws}
From the estimate for $w'$, we see that  on $(0,\infty)$
\[
\frac{3}{r^2}(w-2)(w + 1) < f(u) = \frac{1}{r}w' + \frac{1}{r^2} (w-2)(w + 3)  < \frac{1}{r^2}(w-2)(2w + 3)  < 0.
\]
\end{remark}

\begin{lemma}
The function $\frac{ru'}{u^2} - \frac{2(s_+ - u)}{s_+ u}$ is strictly decreasing and the function $\frac{ru'}{u^3} - \frac{2(s_+^2 - u^2)}{s_+^2 u^2}$ is strictly increasing. In particular, $\frac{2u(s_+^2 - u^2)}{s_+^2\,r} > u' > \frac{2u(s_+ - u)}{s_+\,r}$.
\end{lemma}

\begin{proof} Let us define $\psi = w - \frac{2(s_+ - u)}{s_+} = w - 2 + \frac{2u}{s_+}$. Using the upper bound on $w'$ we have
\begin{align*}
\psi' 
	&= \underbrace{w'}_{=\frac{1}{r}(\chi + w(w-2))} + \frac{2u'}{s_+}
		< \frac{1}{r}w(w - 2) + \frac{2u'}{s_+}
		= \frac{u'}{u}(\psi - \frac{2u}{s_+}) + \frac{2u'}{s_+} 
		= \frac{u'}{u}\psi.
\end{align*}
It follows that $\frac{\psi}{u}$ is decreasing. It is clear that $\psi (\infty)=0$ and therefore $ \frac{\psi}{u} >0$. Since $u>0$ we have 
$$
u' > \frac{2u(s_+ - u)}{r s_+}.
$$
The monotonicity of the other function and upper bound on $u'$ can be proved similarly.
\end{proof}

\begin{lemma}\label{lemma:relfs}
The following inequality holds
\[
\hat f(u) - 3f(u) > -\frac{2}{w}\,f(u).
\]
\end{lemma}

\begin{proof} Let us define
\[
\psi = w'' + \frac{2}{r}w' - \frac{2}{r^2} w\,(w - 2)(w + 3) = (\hat f(u) - 3f(u))w + 2\,f(u).
\]

In the proof we will refer frequently to equations \eqref{eq:A3-1},\eqref{eq:A3-2} and \eqref{eq:A3-3} without explicitly mentioning.
A simple calculation gives
\begin{align*}
\psi'
	&= \Big[\frac{w'}{w} + \frac{2}{r}\Big]\,w'' - \frac{2}{r^2} \frac{w^3 - 5w - 6}{w}\,w'
		  - \frac{4}{r^3}(w - 2)(w + 3) + \frac{b^2}{3}u'\,w.
\end{align*}
Using the following inequality
\[
\frac{b^2}{3}\,r\,u'\,w = \frac{b^2}{3}\,u\,w^2 \leq [\hat f(u) - 3f(u)]\,w^2 = w\,w'' - \frac{2}{r^2}w(w - 2)(w + 1)(w + 3),
\]
we have
\begin{align*}
\psi'
	&\leq \Big[\frac{w'}{w} + \frac{2 + w}{r}\Big]\,\Big(\psi - \frac{2}{r}w' + \frac{2}{r^2}w (w - 2)(w + 3)\Big) - \frac{2}{r^2} \frac{w^3 - 5w - 6}{w}\,w'\\
		&\qquad  - \frac{2}{r^3}(w - 2)(w + 3)(w^2 + w + 2) \\
	&= \Big[\frac{w'}{w} + \frac{2 +w}{r}\Big]\,\psi - \frac{2}{r}\frac{|w'|^2}{w} - \frac{6}{r^2} \frac{w-2}{w}\,w'
		 + \frac{2}{r^3}(w - 2)^2(w + 3)\\
	&< \Big[\frac{w'}{w} + \frac{2 +w}{r}\Big]\,\psi - \frac{2w'}{r^2 w} (w + 3)(w-2)
		 + \frac{2}{r^3}(w - 2)^2(w + 3)\\
	&< \Big[\frac{w'}{w} + \frac{2 +w}{r}\Big]\psi
	,
\end{align*}
where we used $w' < w(w-2)<0$ in the last two estimates. Recalling that $w = \frac{ru'}{u}$, we see that
\[
\psi' <  \Big[\frac{w'}{w} + \frac{2}{r}  + \frac{u'}{u}\Big]\psi, \text{ or equivalently } \frac{d}{dr} \frac{\psi}{r^2\,w\,u} < 0.
\]
It follows that $\frac{\psi}{r^2\,w\,u}$ is a decreasing function. Since  $\frac{\psi}{r^2\,w\,u} \to 0$ as $r \to \infty$ infinity, we conclude that $\psi$ is positive. The statement of the lemma follows.
\end{proof}

As a consequence of the above results, we have the following lower and upper bounds for the solution
\begin{corollary}\label{Cor:adLB}
Assume that ${u(r)} = \alpha{r^2} + o(r^2)$ as $r \to 0$ and $u(r) = s_+ - \beta r^{-2} + o(r^{-2})$ as $r \to \infty$.  Then $u(r)$ has the following upper and lower bounds
\begin{align}
u(r) 
	&\geq \frac{s_+\alpha\,r^2}{\alpha\,r^2 + s_+},\label{Eq:LB1}\\
u(r) 
	&\leq \frac{s_+^2\,r^2}{s_+\,r^2 + \beta},\label{Eq:UB1}\\
u(r)
	&\leq \frac{s_+\alpha\,r^2}{\sqrt{\alpha^2\,r^4 + s_+^2}}.\label{Eq:UB2}
\end{align}
\end{corollary}

\begin{proof} Using Remark~\ref{remark:fws} we have
\[
u'' + \frac{2}{r}u' - \frac{6}{r^2} u = u\,f(u) \leq \frac{2|u'|^2}{u} - \frac{1}{r}u' - \frac{6}{r^2}u.
\]
It follows that $u'' + (-\frac{2u'}{u} + \frac{3}{r})u' \leq 0$, which is equivalent to $\frac{d}{dr} \frac{r^3u'}{u^2} \leq 0$. Integrating this inequality and using the fact that 
$$
\lim_{r \to 0} \frac{r^3 u'}{u^2} =  \frac{2}{\alpha}, \hbox{ and } \lim_{r \to \infty} \frac{r^3 u'}{u^2} =  \frac{2\beta}{s_+^2}
$$
we obtain
\begin{equation}
\frac{2\beta}{s_+^2} \leq \frac{r^3 u'}{u^2} \leq \frac{2}{\alpha}.
	\label{Eq:13VI13.01}
\end{equation}

The second inequality in \eqref{Eq:13VI13.01} implies that $\frac{d}{dr} \Big( \frac{1}{u} - \frac{1}{\alpha r^2}\Big) \geq 0$, and integrating it from $r$ to $\infty$ we obtain $\frac{1}{u} -  \frac{1}{\alpha r^2} \leq \frac{1}{s_+}$,
which implies \eqref{Eq:LB1}. Similarly, the first inequality in \eqref{Eq:13VI13.01} implies that $\frac{d}{dr} \Big( \frac{1}{u} - \frac{\beta}{s_+^2 r^2}\Big) \leq 0$, which leads to \eqref{Eq:UB1}.

To prove \eqref{Eq:UB2}, we again use Remark~\ref{remark:fws}. We have
\[
u'' + \frac{2}{r}u' - \frac{6}{r^2} u = u\,f(u) \geq \frac{3|u'|^2}{u} - \frac{3}{r}u' - \frac{6}{r^2}u.
\]
It follows that $u'' + (-\frac{3u'}{u} + \frac{5}{r})u' \geq 0$, which is equivalent to $\frac{d}{dr} \frac{r^5u'}{u^3} \geq 0$. Using the same argument as before we obtain
\[
\frac{r^5 u'}{u^3} \geq \frac{2}{\alpha^2}, \,\, \text{ i.e. } \,\,\frac{d}{dr}\Big(\frac{1}{u^2} - \frac{1}{\alpha^2\,r^4}\Big) \leq 0.
\]
Consequently
$\frac{1}{u^2} -  \frac{1}{\alpha^2 r^4} \geq \frac{1}{s_+^2}$,
which implies \eqref{Eq:UB2}.
\end{proof}

In Corollary \ref{Cor:adLB}, the lower bound of $u$ depends on $u''(0)$ which is a priori unknown. The following result gives an lower bound which is independent of $u''(0)$.

\begin{lemma} 
There holds
$$u(r) > \underline u(r):= \frac{b^2}{2c^2}\frac{r^6}{(r^2+ \frac{36c^2}{b^4})(r^4 + \frac{12^4 c^4}{b^8})} \text{ for } r \in (0,\infty).$$ 
\end{lemma}

\begin{proof} Let $u_0$ denote the positive solution of \eqref{RS::ODE} corresponding to $a = 0$. Let us observe first that we have:
\begin{equation}\label{u0a}
u_0(r)\le u(r), \forall r>0.
\end{equation}
This follows from the comparison principle in Proposition~\ref{Prop:CompPrin}, Remark \ref{RemW}, the fact that $u_0$ is a sub-solution of \eqref{RS::ODE} for $a > 0$ (in the sense of Definition~\ref{def1}), and
\[
u(\infty) = s_+ = \frac{b^2 + \sqrt{b^4 + 24a^2c^2}}{4c^2} > \frac{b^2}{2c^2} = u_0(\infty) \text{ whenever } a > 0.
\]

Therefore, it suffices to show that $u_0 \geq \underline{u}$. Using Remark~\ref{rmk:physrescaling} in the introduction, it suffices to check this for e.g. $b = c = 1$. In that case, a lengthy computation shows that
\begin{multline*}
\underline{u}'' + \frac{2}{r}\underline{u}' - \frac{6}{r^2}\underline{u} - f_0(\underline{u})\underline{u} = \frac{36r^4}{(r^2+36)^3(r^4+12^4)^3} \times \\
	\times (278628139008 + 9029615616 r^2+ 85100544 r^4-373248 r^6-5184 r^8+41 r^{10}) , 
\end{multline*}
where $f_0(u) = - \frac{1}{3}u + \frac{2}{3}u^2$. It is straightforward to check that 
$$
278628139008 + 9029615616 r^2+ 85100544 r^4-373248 r^6-5184 r^8+41 r^{10} >0 \hbox{ on } (0,\infty).
$$
In other words the function $\underline u$ is a sub-solution of \eqref{RS::ODE} with $a = 0$.  

Notice that
$$
\underline u (r) = \frac{1}{2} - \frac{18}{r^2} + o(r^{-2}) \text{ and } u_0(r) = \frac{1}{2} - \frac{18}{r^2} + o(r^{-2})   \text{ as } r\rightarrow \infty,
$$
where we have used Proposition \ref{ERS::Asymp}. Taking into account the behaviour of $\underline u$ at $0$ we can apply again Proposition~\ref{Prop:CompPrin} to obtain that $\underline u(r)\le u_0(r),\forall r>0$.
\end{proof}

\section{Existence of sign-changing solutions}
\label{sec:SC}

In Section \ref{sec:uniqueness}, we show that, for $F$ satisfying \eqref{condF}, the problem \eqref{RS::ODE}\&\eqref{BC} has a unique positive solution. Furthermore, under more stringent conditions on $F$, that solution is the unique solution of the problem \eqref{RS::ODE}\&\eqref{BC}. The goal of this section is to give examples of nonlinearities $F$ (which satisfy \eqref{condF}) such that, for any finite interval $(0,R)$, the problem \eqref{RS::ODE}\&\eqref{BC} has another solution besides the positive solution. This additional solution is necessarily sign-changing (in view of Theorem \ref{thm:main}) and is of mountain-pass type.

For simplicity, we set
\[
p = 2 \, \, \text{ and } \, \, q = 6.
\]
The problem \eqref{RS::ODE}\&\eqref{BC} becomes
\begin{align*}
&u'' + \frac{2}{r} u' -\frac{6}{r^2} u = F(u) \text{ in } (0,R),\\
&u(0) = 0, u(R) = s_+.
\end{align*}
Let $u_*$ be the positive solution obtained in Theorem \ref{thm:main}.

\subsection{Minimizing properties}

We have seen in Corollary~\ref{cor:energyuniqueness}, proved in Section~\ref{ssec:NoPosAs} that if $F$ satisfies \eqref{Fcond:even},
then for $R\in(0,\infty)$ the function $u_R$ (the solution of  \eqref{RS::ODE}\&\eqref{BC} obtained in Theorem~\ref{thm:main}) is actually a global minimizer of the energy $E$ defined in \eqref{energyODE}, in the introduction. It is natural to ask if $u_R$ is actually a global minimizer for $E$ when $F$ does not necessarily satisfy \eqref{Fcond:even}. In general the answer is negative. For example, for the nonlinearity $F(u) = u^4 - u$, the energy $E$ is unbounded from below. However, we prove:

\begin{lemma}\label{lem:locminimality}
Assume $p = 2$, $q = 6$, $R \in (0,\infty)$ and $F$ satisfies \eqref{condF}. Let $u_*$ be the positive solution in Theorem \ref{thm:main}. Then $u_*$ is a strictly stable local minimizer for $E[\cdot;(0,R)]$.
\end{lemma}

\begin{proof} Consider the second variation of $E$ at $u_*$:
\[
Q[v] := \int_0^R \Big[r^2|v'|^2 + 6\,v^2 + r^2\,F'(u_*)\,v^2\Big]\,dr,
\]
where $v$ belongs to
\[
\mcM_0 := \Big\{v : (0,R)\rightarrow \RR \,\Big| rv' \in L^2(0,R), v(R) = 0\Big\}.
\]

It suffices to prove that, for some $\delta > 0$,
\[
Q[v] \geq \delta \int_0^R r^2 |v'|^2\,dr =: \delta\|v\|^2 \text{ for all } v \in \mcM_0.
\]
By a standard density argument, it suffices to prove the above for $v \in C_c^\infty(0,R)$.

To this end, note that $u_*'$ satisfies
\[
(u_*')'' + \frac{2}{r}(u_*')' - \frac{8}{r^2}u_*' + \frac{12}{r^3}u_* = F'(u_*)u_*' \text{ in } (0,R).
\]
Furthermore, by Lemma \ref{ERS:monotonicity}, $u_*'$ is non-negative.

Fix $v \in C_c^\infty(0,R)$ and write $v = u_*'\,w$. We have
\begin{align*}
Q[v] 
	&= \int_0^R \Big[r^2|(u_*' w)'|^2 + 6|u_*'|^2 w^2 + r^2 F'(u_*)|u_*'|^2\,w^2\Big]\,dr\\
	&= \int_0^R \Big[r^2|(u_*' w)'|^2 + 6|u_*'|^2 w^2 + ((r^2u_*'')' - 8 u_*' + \frac{12}{r}u_*)u_*'\,w^2\Big]\,dr\\
	&= \int_0^R \Big[r^2|u_*'|^2 |w'|^2 - 2|u_*'|^2 w^2 + \frac{12}{r}u_*\,u_*'\,w^2\Big]\,dr.
\end{align*}
Recalling \eqref{bound:derivgeneral} and noting that $\gamma_+ = 2$, we obtain
\[
Q[v] 
	\geq \int_0^R \Big[r^2|u_*'|^2 |w'|^2 + 4|u_*'|^2 w^2\Big]\,dr \geq \int_0^R 4v^2.
\]

Since $0 \leq u_* \leq s_+$, $F'(u_*) \geq - C_0$ for some $C_0$ depending only on $F$. It thus follows that
\[
\|v\|^2 \leq Q[v] - \int_0^R r^2\,F'(u_*)\,v^2\,dr \leq Q[v] + C_0R^2\,\int_0^R v^2 \leq \frac{1}{4}(4 + C_0R^2) Q[v],
\]
as desired.
\end{proof}

\subsection{Mountain pass solutions}

In this subsection we obtain a mountain-pass solution for the BVP \eqref{RS::ODE}\&\eqref{BC} on finite domains when the nonlinearity $F$ satisfies a certain growth condition.

\begin{proposition}\label{lem:MP}
Assume $p = 2$, $q = 6$, $R \in (0,\infty)$. Assume that $F$ satisfies \eqref{condF} and, for some $\varkappa > 0$, $0 \leq \lambda < 4$ and $C > 0$ we have $F(t) = \varkappa t^4 + \mathring{F}(t)$ with $\mathring{F}$ satisfying
\begin{equation}
|\mathring{F}(t)|
	 \leq C(1 + |t|^\lambda) \text{ for } t \in \RR,
	\label{FCond:quartic}
\end{equation}

Then besides the positive solution obtained in Theorem \ref{thm:main}, the problem \eqref{RS::ODE}\&\eqref{BC} admits a sign-changing solution.
\end{proposition}

For example, we note that the nonlinearities $F(u) = u^4 + 2u^3 - u^2 - 2u$ and $F(u) = u^4 - u^3$ satisfy all hypotheses of Lemma \ref{lem:MP}.

\begin{proof}  Let us consider the set
\[
\mcM:=\Big\{u:(0,R)\to\RR\, :\, ru', u\in L^2(0,R),\quad u(R)=s_+\Big\}.
\]

It is easy to check that $u \in \mcM$ is a critical point for $E$ if and only if $v = u_* - u \in \mcM_0$ is a critical point of
\[
I[v] = \frac{1}{2}\int_0^R \Big[r^2 |v'|^2 + 6v^2 + r^2(2F(u_*)v + h(u_* - v) - h(u_*))\Big]\,dr,
\]
where $h$ is given by \eqref{Eq:hDefinition}.

Note that $\mcM_0$ is a Hilbert space with respect to the inner product
$\langle v_1, v_2\rangle = \int_0^R r^2\,v_1'\,v_2'\,dr$. The (Fr\'echet) derivative of $I$ is given by
\begin{equation}
\langle I'[v],\varphi\rangle = \int_0^R  \Big[r^2 v' \varphi' + 6v \varphi + r^2(F(u_*) - F(u_* - v))\varphi \Big]\,dr.
	\label{Eq:IprimeDef}
\end{equation}

By Lemma \ref{lem:locminimality}, $0$ is a strictly stable local minimizer of $I$ and $I[0] = 0$. In addition, for $v \geq 0$ and $v \not\equiv 0$, we have $I[tv] \rightarrow -\infty$ as $t \rightarrow \infty$ thanks to \eqref{FCond:quartic}.  We would like to find a second critical point of $I$ via the mountain pass theorem (see e.g. \cite{Rabinowitz}). To this end, it remains to show that $I$ satisfies the Palais-Smale condition. More precisely, we need to show that, if $v_n$ is a sequence in $\mcM_0$ satisfying
\[
I[v_n] \leq C \text{ and } I'[v_n] \rightarrow 0,
\]
then $v_n$ has a convergent subsequence. Note that, by standard elliptic estimates, it suffices to show that the sequence $v_n$ is bounded in $\mcM_0$. 

Let $A_n = \int_0^R [r^2 |v_n'|^2 + 6v_n^2]\,dr$ and fix some $\delta > 0$ small.

First, taking $\varphi = v_n^+$ in \eqref{Eq:IprimeDef} and noting that $I'[v_n] \rightarrow 0$, we can find some $\epsilon_n \rightarrow 0$ such that
\[
-\epsilon_n \sqrt{A_n}
	 \leq  \int_{\{v_n > 0\}}  \Big[r^2 |v_n'|^2 + 6v_n^2 + r^2(F(u_*) - F(u_* - v_n))v_n \Big]\,dr.
\]
Thus, by \eqref{FCond:quartic},
\begin{equation}
-\epsilon_n \sqrt{A_n}
	 \leq  \int_{\{v_n > 0\}}  \Big[r^2 |v_n'|^2 + 6v_n^2 - (1-\delta) \varkappa\, r^2 |v_n|^5 \Big]\,dr + C,
	 \label{MP:PositivePart}
\end{equation}
where here and below $C$ denotes some constant that may vary from line to line but is always independent of the sequence $v_n$.

Next, using the boundedness of $I[v_n]$ and \eqref{FCond:quartic},
\begin{align*}
C
	&\geq A_n + \int_0^R \big[2F(u_*)v_n + h(u_* - v_n) - h(u_*)\big]\,r^2\,dr\\
	&\geq A_n - \int_{\{v_n > 0\}} \frac{2\varkappa}{5}(1 + \delta)\,r^2\,v_n^5\,dr - C.
\end{align*}
Thus, by \eqref{MP:PositivePart},
\[
C \geq \frac{3}{5}(1 - O(\delta))A_n - \frac{2}{5}(1 + O(\delta))\epsilon_n\,\sqrt{A_n}.
\]
This implies the boundedness of $A_n$ as desired. The mountain pass theorem can then be invoked to assert the existence of a second critical point of $I$, thus of $E$, which is a solution of \eqref{RS::ODE}\&\eqref{BC}. Since positive solution of \eqref{RS::ODE}\&\eqref{BC} is unique, this second solution must be sign-changing.
\end{proof}

\appendix

\section{Lifting for radially symmetric $Q$-tensors}
\label{sect:ode_pde}

In this appendix, we classify radially symmetric matrix-valued maps by using only one degree of freedom, the scalar $u(|x|)$, also called lifting.

\begin{lemma}\label{RSForm}
If $Q: B_R(0) \rightarrow \mcS_0$ is a radially symmetric measurable map, then there exists a measurable function $u: (0, R) \rightarrow \RR$ such that
\begin{equation}
Q(x) = u(|x|)\,\hh (x)
	\quad \text{ for a.e. }\,  x \in B_R(0)\;.
\label{RSF::Rep}
\end{equation} where $\hh (x):=\left(\frac{x}{|x|}\otimes \frac{x}{|x|}-\frac{1}{3}Id\right)$.
The function $u$ is given by 
\be
\label{ident_u_rad}u(|x|) = {3 \over 2} \tr (Q(x) \hh (x)) \quad \textrm{ a.e. in } \quad B_R(0).\ee If the origin $0$ is a Lebesgue point of $Q$, then it is also a Lebesgue point of $u$, and $Q(0)=0$, $u(0)=0$. Moreover, if $Q$ is continuous on $B_R(0)$, then $u$ is also a continuous function on $[0, R)$ with $u(0)=0$.
\end{lemma}

\begin{proof} Fix a point $x \in B_R(0)$ where \eqref{RadSymDef} holds. Write $x=rp$ for some $p \in \Sphere^2$ and $r \geq 0$. Assume for now $x \neq 0$. Let $G_x$ denote the subgroup of rotation matrices in $SO(3)$ that fixes $x$, i.e. $\Rot x =x$ for all $R \in G_x$. By the definition of radial symmetry for tensors, we have
\begin{equation} \label{Q-sym}
Q(x) = Q(\Rot x) = \Rot\,Q(x)\,\Rot^t \text{ for any } \Rot \in G_x
	\;.
\end{equation}
Observe that for $x\neq 0$, $\Spec (\hh(x))=\{-\frac 1 3, -\frac 1 3, \frac 2 3\}$ and the eigenspaces corresponding to the eigenvalues $-1/3$ and $2/3$ of $\hh(x)$ are given by the plane $(\RR x)^\perp$ and the line $\RR x$, respectively. In view of \eqref{RSF::Rep}, it is then natural to prove that:
\begin{claim}\nonumber
\label{claim_eigen}
 $Q(x)$ cannot have three distinct eigenvalues. Moreover, $Q(x) = u(x) \hh (x)$ for some $u(x) \in \RR$.
 \end{claim}
 \begin{proof}
First, observe that 
\begin{multline}
\text{if $v$ is an eigenvector of $Q(x)$ then $\Rot v$ is also an eigenvector of $Q(x)$}\\
	\text{for all $\Rot \in G_x$ (with the same eigenvalue).}
	\label{RotInvProp}
\end{multline}
Indeed, $Q(x) v =\lambda v$ in view of \eqref{Q-sym} implies
$$
Q(x) \Rot v = \Rot\,Q(x)\,\Rot^t \Rot v =\lambda \Rot v.
$$
To prove our Claim, we distinguish the following two cases (since $Q(x)$ is a symmetric matrix, so that $\RR^3$ is a direct sum of eigenspaces of $Q(x)$):

\medskip

\noindent{\it Case 1: $Q(x)$ has an eigenvector $v$ which is neither parallel nor perpendicular to $x$.} Then \eqref{RotInvProp} implies that the whole $\bR^3$ is an eigenspace of $Q(x)$ corresponding to a single eigenvalue. Since $Q(x)$ is traceless we deduce that $Q(x)=0$, i.e., all eigenvalues of $Q(x)$ are zero.

\medskip

\noindent {\it Case 2: $Q(x)$ has an eigenvector $v$ which is parallel to $x$ and two linear independent eigenvectors $v_2$ and $v_3$ which are perpendicular to $x$.} Let $\lambda_1$, $\lambda_2$ and $\lambda_3$ be the corresponding eigenvalues. Then \eqref{RotInvProp} implies that $\lambda_2 = \lambda_3$. By tracelessness of $Q(x)$, $\lambda_2 = \lambda_3 = -\frac{1}{2}\lambda_1$. Furthermore, $Q(x)$ has the same eigenspaces as $ \hh (x)$ so that $Q(x) = u(x) \hh (x)$ for some $u(x) \in \RR$ (here, $u(x)=-3\lambda_2$).

\medskip

\noindent In both cases, we obtain the representation $Q(x) = u(x) \hh (x)$ which proves our Claim.
\end{proof}

\noindent To finish the proof of our lemma, notice that $ \hh (\Rot x)=\Rot \hh (x) \Rot^t$ for all $\Rot \in SO(3)$ and \eqref{RadSymDef} also holds at every point $\tilde x=\tilde \Rot x$ for every rotation $\tilde \Rot \in SO(3)$ (since \eqref{RadSymDef} holds at $x$ by our assumption). Combined with Claim \ref{claim_eigen}, it follows that $u(x)=u(\Rot x)$ for all $\Rot \in SO(3)$ which entails that $u$ is indeed a function of $|x|$, i.e., \eqref{RSF::Rep} holds at $x$.
From here, it is easy to see that \eqref{ident_u_rad} also holds at $x$.

Assume now that $x = 0$ is a Lebesgue point of $Q$, i.e., there exists a matrix $Q^*$ such that 
$$\lim_{r\to 0} -\hspace{-4.3mm} \int_{B_r(0)} |Q(x)-Q^*|\, dx=\lim_{r\to 0} -\hspace{-4.3mm} \int_{B_r(0)} |Q(\Rot x)-Q^*|\, dx=0$$ by the change of variable $\tilde x:=\Rot x$ for some $\Rot \in SO(3)$. Since for a.e. $x\in B_R(0)$, $Q(x)\in \mcS_0$ we deduce that $Q^* \in \mcS_0$.
Since \eqref{RadSymDef} holds a.e. in $B_R(0)$, we deduce that 
$$\lim_{r\to 0} -\hspace{-4.3mm} \int_{B_r(0)} |Q(x)-\Rot^tQ^*\Rot|\, dx=0,$$ so that
$Q^*=\Rot^t Q^*\Rot$ for all $\Rot \in SO(3)$. Since $Q^*$ is a traceless symmetric matrix, it implies that $Q^*=0$. Relation \eqref{ident_u_rad} allows to obtain that $0$ is also a Lebesgue point for $u$ and $u(0)=0$. For the last assertion, assume that $Q$ is continuous. Obviously, \eqref{RadSymDef} holds everywhere in $B_R(0)$. By  \eqref{ident_u_rad}, since $\hh$ is continuous away from $0$ and bounded near $0$, the continuity of $u$ on $(0, R)$ immediately follows. Since $Q$ is assumed to be continuous at $0$, by \eqref{ident_u_rad}, we deduce that $u$ can be continuously extended to $u:[0, R)\to \RR$ by setting $u(0)=0$. 
\end{proof}

\section{Some maximum principles}
\label{app}

In this appendix, we present some maximum principles which were needed in the body of the paper.

\begin{lemma}\label{lemma:weakmax}
For $R\in (0, \infty]$, $p,q\in C(0,R)$ with $q(r)\ge 0,\forall r\in (0,R)$ we denote 
$$Lw:=-w''-p(r)w'+q(r)w.$$ 
Assume that there exists a nonnegative function $w_0\in C^2(0,R)$ with
$Lw_0\geq0$ and $\lim_{r\to 0}w_0(r)=\infty$. If $w\in W_{loc}^{1,\infty}(0,R)\cap L^\infty(0,R)$  and  $Lw\ge 0$ in the sense of distributions (or measures) in $(0,R)$ with $\liminf_{r \rightarrow R} w(r)\ge 0$ then $$w(r)\ge 0,\forall r\in (0,R).$$
\end{lemma}
\begin{proof} The result of the lemma and its proof are well-known to  experts, but we provide the proof here for completeness. We pick an arbitrary $\varepsilon>0$. There exists $\delta_0(\varepsilon), M_0(\eps)\in (0, R)$ so that 
\be
\label{cond_bdry3}
w(\delta)\geq-\varepsilon w_0(\delta)> -\varepsilon\,w_0(\delta) - \varepsilon \quad \textrm{ and } \quad w(M) >-\eps\geq -\varepsilon\,w_0(M) - \varepsilon\ee for all $0 < \delta<\delta_0(\varepsilon)$ and $M_0(\eps)<M<R$. Hence, by the usual weak maximum principle applied to $L$ on the interval $(\delta, M)$ we get $w\ge -\varepsilon w_0 - \varepsilon$ in $(\delta, M)$ for any $0 < \delta < \delta_0(\eps)$ and $M_0(\eps)<M<R$. Indeed, if we denote by $v:=w+\varepsilon w_0 + \varepsilon$ we obtain that $Lv\geq 0$ in the sense of distributions in $(0,R)$. Set $v_-:=\max\{0, -v\}$ and $P\in C^1(0, R)$ be a primitive of $p$, i.e., $P'=p$ on $(0, R)$. Noting that $0\leq e^{P} v_-\in C_c\big((\delta, M)\big)$ (due to \eqref{cond_bdry3}), and using it as a test function we obtain
$$0\geq \int_\delta^{M} \bigg((v_-')^2+q(r)(v_-)^2\bigg)e^{P(r)}\, dr$$
and conclude that $v_-\equiv 0$ on $(\delta, M)$.
 Since we can choose any $\delta\in (0,\delta_0(\varepsilon))$ and $M\in (M_0(\eps), R)$, we have in fact $w\ge -\varepsilon w_0 - \varepsilon$ in $(0,R)$. Since $\varepsilon>0$ was arbitrary we can let $\varepsilon\to 0$ and obtain the conclusion.
\end{proof}

\begin{remark}\label{remark:fundsol}
The above maximum principle was used in two specific cases, namely:
\begin{itemize}
\item For $p(r)=\frac{p}{r},q(r)=\frac{q}{r^2}$ with $w_0(r)=r^{\gamma_-}$ where $\gamma_-$ is the negative Fuchsian index of (\ref{RS::ODE}) (see \eqref{Fuch_ind}).
\item For $p(r)=-\frac{p-2}{r},q(r)=\frac{A}{r^4}$, $A>0$ with $w_0(r)=r^B I_B(\frac{\sqrt{A}}{r})$ where $B=(p-1)/2$ and $I_B$ is the modified Bessel function (see for instance \cite{Abramowitz}, p. $375$, $9.6.10$ and p. $377$, $9.7.1$) that satisfies the modified Bessel's equation
$$I_B''(t)+\frac 1 t I_B'(t)-(1+\frac{B^2}{t^2})I_B(t)=0, \quad t>0.$$  
and  has exponentially growth at infinity.
\end{itemize}
\end{remark}

\begin{lemma}\label{Lem:SingHopf}
Assume that $w \in C^2(0,R) \cap C[0,R)$ satisfies
\[
Lw(r) := w''(r) + \frac{a}{r}\,w'(r) + c(r)\,w(r) \leq 0 \text{ in } (0,R)
\]
for some constant $a \geq 1$ and some function $c \in C[0,R)$. If $w \geq 0$ in $(0,R)$ and $w(0) = 0$, then $w \equiv 0$.
\end{lemma}

\begin{remark}
The conclusion is not true for $0 < a < 1$. For example, take $w(r) = r^{1-a}$ and $c \equiv 0$.
\end{remark} 

\begin{proof} Arguing by contradiction, assume that $w \not\equiv 0$. By the standard strong maximum principle, $w > 0$ in $(0,R)$.

For some small positive $\epsilon$, consider the function 
\[
\psi(r) = \epsilon^2 - (\epsilon - r)^2.
\]
We have, for $0 < \delta < \epsilon$,
\[
L(\psi - \delta) = -2 + \frac{2a(\epsilon - r)}{r} + c(r)[\epsilon^2 - \delta - (\epsilon - r)^2].
\]
Since $a > 0$, there exists some $\epsilon>\lambda > 0$ independent of $\delta$ such that $L\psi > 0$ in $(0,\lambda)$. 

Pick $\mu > 0$ such that $\mu\psi(\lambda) < w(\lambda)$. Clearly there exists some $0 < \lambda' < \lambda$ which might depend on $\delta$ such that $w \geq \mu(\psi - \delta)$ in $[0,\lambda']$. By the maximum principle, $w \geq \mu (\psi-\delta)$ in $[\lambda',\lambda]$. It follows that $w \geq \mu(\psi - \delta)$ in $[0,\lambda]$ for all $0 < \delta < \epsilon$, which implies that $w \geq \mu\psi$ in $[0,\lambda]$. Since $w(0) = \psi(0) = 0$, this implies that
\[
\liminf_{r \rightarrow 0} \frac{w(r)}{r} \geq \mu\psi'(0) = 2\mu\epsilon > 0.
\]
In particular, there exists a sequence $r_k \rightarrow 0$ such that
\[
w'(r_k) \geq \mu\epsilon > 0.
\]

Recall that $w$ and $c$ are continuous up to $r = 0$ and $w(0) = 0$. Thus, we can choose some $\eta > 0$ such that $|c(r)w(r)| < \frac{(a+1)\mu\epsilon}{2\eta}$ for $0 < r < \eta$. Thus, as $Lw \leq 0$,
\[
(r^a w')' \leq \frac{(a+1)\mu\epsilon}{2\eta} r^a \text{ for } 0 < r < \eta.
\]
Fix some $r_k < \eta$. Then
\[
w'(r) \geq \frac{1}{r^a}\Big[r_k^a\,w'(r_k) - \frac{\mu\epsilon}{2\eta}r_k^{a+1}\Big] \geq \frac{\mu\epsilon}{2} \frac{r_k^a}{r^a} \text{ for } 0  < r < r_k.
\]
Since $a \geq 1$, this implies that
\[
\lim_{r \rightarrow 0} w(r) = - \infty,
\]
contradicting our hypothesis that $w$ is continuous up to $r = 0$ and $w(0) = 0$.
\end{proof}

\section*{Acknowledgments} The authors gratefully acknowledge the hospitality of Hausdorff Research Institute for Mathematics, Bonn where part of this work was carried out.
R.I.  acknowledges partial support by the ANR project ANR-10-JCJC 0106. V.S. acknowledges support by EPSRC grant  EP/I028714/1.

\bibliographystyle{siam}
\bibliography{paris,%
LiquidCrystals}

\end{document}